\documentclass[11pt]{article}
\usepackage{amsmath}
\usepackage{amssymb}
\usepackage{amsfonts}
\usepackage{mathrsfs}
\usepackage{shortvrb,psfrag}
\usepackage{enumerate}
\usepackage{color}
\usepackage[dvips]{graphicx}

\let\f=\frac

\def\ba{\begin{eqnarray}}
\def\ea{\end{eqnarray}}

\def\mne{\mn_{\epsilon}}
\def\je{{\cal J}_{\epsilon}}
\def\mn{{\mathbf{n}}}

\def\mcurl{{\mathbf{ curl}}}

\def\mh{{\mathbf{h}}}

\def\cB{{\cal B}}

\def\cM{{\cal M}}

\def\R{\Bbb R}

\def\no{\noindent}
\def\na{\nabla}

\def\eqdef{\buildrel\hbox{\footnotesize def}\over =}
\def\endproof{\hphantom{MM}\hfill\llap{$\square$}\goodbreak}

\newcommand{\beq}{\begin{equation}}
\newcommand{\eeq}{\end{equation}}
\newcommand{\ben}{\begin{eqnarray}}
\newcommand{\een}{\end{eqnarray}}
\newcommand{\beno}{\begin{eqnarray*}}
\newcommand{\eeno}{\end{eqnarray*}}

\newtheorem{Theorem}{Theorem}[section]

\newtheorem{Proposition}[Theorem]{Proposition}
\newtheorem{Lemma}[Theorem]{Lemma}
\newtheorem{Corollary}[Theorem]{Corollary}
\newtheorem{Remark}[Theorem]{Remark}

\begin{document}

\title{The Landau-Lifshitz equation
of the ferromagnetic spin chain
and Oseen-Frank flow}

\author{ Xueke PU, Meng WANG, Wendong WANG}


\date{\today}
\maketitle

\begin{abstract}

\textcolor[rgb]{0.00,0.00,0.00}{In this paper, we consider the Landau-Lifshitz equation of the ferromagnetic spin chain from $\R^2$ to the unit sphere $S^2$ under the general Oseen-Frank energy. We obtain global existence and uniqueness of weak solutions for large energy data; moreover, the number of singular points is finite.}
\end{abstract}

\setcounter{equation}{0}

\section{Introduction}
The \textcolor[rgb]{0.00,0.00,0.00}{$d$-dimensional} classical system for the isotropic Heisenberg chain with spin vector $\mn=(n_1,n_2,n_3)$ is described by the Hamiltonian density (without external magnetic field) $H=|\nabla\mn|^2/2$. The spin equation of motion with the Gilbert damping term (without the external magnetic field) has the form
\begin{equation}\label{p1}
\partial_t\mn= \alpha\mn\times(\mn\times \frac{\delta H}{\delta\mn}) -{\beta}\mn\times \frac{\delta H}{\delta\mn},
\end{equation}
where $\alpha\geq0$ is the Gilbert damping constant and $\beta$ is the exchange constant satisfying $\alpha^2+\beta^2=1$ and $H$ is the Hamiltonian density. Explicitly, this gives the following classical Landau-Lifshitz equation
\ben\label{eq:classical Landau-Lifshitz}
\partial_t\mn=\beta\mn\times \Delta\mn -\alpha\mn\times (\mn\times\Delta\mn).
\een
The above system \eqref{p1} or \eqref{eq:classical Landau-Lifshitz} is called the Landau-Lifshitz equation or the Landau-Lifshitz-Gilbert equation, which was first derived on phenomenological grounds by Landau-Lifshitz in \cite{LandauL1935}. It gives rise to a continuum spin wave theory. Note that the above system (\ref{eq:classical Landau-Lifshitz}) reduces to the heat flow of harmonic maps when $\alpha=1,\beta=0$ and to the Schr\"{o}dinger flow when $\alpha=0,\beta=1$.


Motivated by the study on the heat flow of harmonic maps (see \cite{ES1964,SchoenU82,Struwe85,Struwe88,Chenstruwe1989} and so on) and Schr\"{o}inger flow (see \cite{DingWang1998,DingWang2001,TerngUhlen2006} and so on), much progress has been made recently in the analysis of the Landau-Lifshitz-Gilbert Equation \eqref{eq:classical Landau-Lifshitz}. For example, see \cite{Alouges92} for the existence of global weak solutions of \eqref{eq:classical Landau-Lifshitz} under the Neumann boundary condition in any dimensions, and see \cite{CarbouF2001,Moser1,DingGuo2004,Liu2004,Moser2,Atsushi2005,DingWang2007,GuoDing2008}  and the references therein for partial regularity and the analysis of singularity of the system \eqref{eq:classical Landau-Lifshitz}. More recently, the existence of partially smooth, global weak solutions of \eqref{eq:classical Landau-Lifshitz} similar to \cite{Struwe85}, has been obtained by Guo-Hong \cite{GuoHong1993} for $d=2$, Melcher \cite{Melcher2005} for $d=3$, and Wang \cite{Wang2006} for $d=4$ with Dirichlet boundary conditions.  More recently, the first author and Guo \cite{PG11,PG13} studied the fractional generalization of the Landau-Lifshitz equation and obtained local well-posedness and global existence of weak solutions.


\textcolor[rgb]{0.00,0.00,0.00}{In this paper, we shall consider the case when the energy density is replaced by the Oseen-Frank energy density. The Oseen-Frank energy density expresses the free energy density of a nematic liquid crystal in terms of its optic axis, and is a measure of the increase in the Helmholtz free energy per unit volume due to deviations in the orientational ordering away from a uniformly aligned nematic director configuration. See \cite{HartKL1986} for the analysis for the minimizers of the Oseen-Frank energy. Let $W=W(\mn,\nabla\mn)$ be the Oseen-Frank density of the form
\begin{align*}
W(\mn,\nabla\mn)=&k_1({\rm div}\mn)^2+k_2|\mn\times({\mathbf{\nabla}}\times\mn)|^2+k_3|\mn\cdot({\mathbf{\nabla}}\times\mn)|^2 \\
&+(k_2+k_4)\big(\textrm{tr}(\nabla\mn)^2-({\rm div}\mn)^2\big),
\end{align*}
where $k_1,k_2,k_3,k_4$ are elastic constants depending on the materials and temperature.}

Replacing $H$ in \eqref{p1} with $W$, we obtain the Landau-Lifshitz equation of Oseen-Frank energy as follows:
\begin{equation}\label{eq: Landau-Lifshitz with Oseen-Frank}
\partial_t\mn= -\alpha\mn\times(\mn\times\mh) +\beta\mn\times\mh,
\end{equation}
where the vector field $\mh$ is given by
$$\mh=-\frac{\delta W}{\delta \mn}=(\nabla_{i}W_{{\rm p}_{i}^{l}}-W_{n_{l}}),$$
where ${\rm p}_{i}^{l}=\nabla_i n_l$ and we adopt the standard summation convention.
Throughout this paper, we denote
\beno
W_{n_i}=\frac {\partial W(\mn,{\rm p})} {\partial {n_i}},\quad W_{{\rm p}_i^l}=\frac {\partial W(\mn, {\rm p})} {\partial {{\rm p}_i^l}}.
\eeno

In what follows, we give explicit form of the vector field $\mh$. For this, we rewrite $W(\mn,\nabla\mn)$ as in \cite{HartKL1986}
\beno
W(\mn,\nabla\mn)=a|\nabla\mn|^2+V(\mn,\nabla\mn),
\eeno
where $a=\min\{k_1,k_2,k_3\}$ and
$$
V(\mn,\nabla\mn)=(k_1-a)({\rm div}\mn)^2+(k_2-a)|\mn\times({\mathbf{\nabla}}\times\mn)|^2
+(k_3-a)|\mn\cdot(\mathbf{\nabla}\times\mn)|^2.
$$
In this way, we have the following (see \cite{WangW2014})

\begin{Lemma}\label{lem:h decomposition}
It holds that
\begin{align*}
(\nabla_{\alpha}W_{p_{\alpha}^{l}})=&2a\Delta\mn+2(k_1-a)\nabla{\rm div}\mn
-2(k_2-a){\mcurl}{(\mn\times({\mcurl}\mn\times\mn))}\nonumber\\
&-2(k_3-a){\mcurl}({\mcurl}\mn\cdot \mn\mn),\nonumber\\
(W_{n_{l}})=&2(k_3-k_2)({\mcurl}\mn\cdot\mn)({\mcurl\mn}),
\end{align*}
In particular, we have
\begin{align}\label{p2}
\mh=&2a\Delta\mn+2(k_1-a)\nabla{\rm div}\mn -2(k_2-a){\mcurl}({\mcurl\mn})\nonumber\\
&-2(k_3-k_2){\mcurl}({\mcurl}\mn\cdot \mn\mn) -2(k_3-k_2)({\mcurl}\mn\cdot\mn)({\mcurl\mn}).
\end{align}
\end{Lemma}

In particular, when $k_1=k_2=k_3$, \eqref{eq: Landau-Lifshitz with Oseen-Frank} with \eqref{p2} reduces to the classical Landau-Lifshitz equation \eqref{eq:classical Landau-Lifshitz}.

In \cite{HongXin2012}, Hong-Xin proved that global existence of weak solution for the  Oseen-Frank flow in 2D (i.e. $\alpha=1,\beta=0$ in (\ref{eq: Landau-Lifshitz with Oseen-Frank})) whose singular points are finite and the uniqueness of weak solution was obtained by the later two authors of \textcolor[rgb]{0.00,0.00,0.00}{the present paper} in \cite{WangWZ2013} (see also \cite{LiTitiXin2014} for different assumptions).

We are aimed to generalize the above results to the general Landau-Lifshitz equation  (\ref{eq: Landau-Lifshitz with Oseen-Frank}) with $\alpha,\beta>0$.
Note that $\partial_{x_3}\mn=0$ in the 2-D case. Let $b\in S^2$ be a constant vector and we define
$$
H_{b}^{1}(\R^2;S^2)=\big\{u: u-b\in H^{1}(\R^2;\R^3), |u|=1 \,\,{\rm a.e. }\,\,\mbox{in} \,\, \R^2 \big\}.
$$
our main results state as follows.
\begin{Theorem}\label{thm:global}
Assume that  the initial data $\mn_0\in H_{b}^{1}(\R^2;S^2)$. Then there exists a unique global weak solution $\mn$ of the system \eqref{eq: Landau-Lifshitz with Oseen-Frank}, which is smooth in $\R^2\times((0,+\infty)\setminus \{T_{i}\}_{i=1}^{L})$ with a finite number of singular points $(x_{i}^{l},T_{i})$, $1\leq l\leq L_i$. Moreover, there are two constants $\epsilon_0>0$ and $R_0>0$ such that each singular point $(x_{i}^{l},T_{i})$ is characterized by the condition
$$
\limsup_{t\uparrow T_{i}}\int_{B_{R}(x_{i}^{l})}|\nabla\mn|^2(\cdot,t)dx> \epsilon_0
$$
for any $R>0$ with $R\le R_0$.
\end{Theorem}

\begin{Remark}
The above theorem generalizes the existence and uniqueness results of the equation (\ref{eq:classical Landau-Lifshitz}) in \cite{GuoHong1993}, and also generalize the existence result in \cite{HongXin2012}. The main difference is the introduced Oseen-Frank energy, which makes the system \eqref{eq: Landau-Lifshitz with Oseen-Frank} does not keep the parabolic property. By constructing strong solutions of a new approximate system, we obtain the local well-posedness and global weak solutions of \eqref{eq: Landau-Lifshitz with Oseen-Frank}. Different with \cite{Struwe85,GuoHong1993}, it's not easy to obtain the uniqueness as said in \cite{HongXin2012}, since the positivity of the diffusion term $\delta_{\mh}\times\mn$ under the metric of $L^2$ norm is unknown. Instead, we
introduce a type of weak Oseen-Frank  metric as in \cite{WangWZ2013}. Our goal is to combine the work of Oseen-Frank energy and the Schr\"{o}nger part $\mn\times\mh$ together.
\end{Remark}


The rest of the paper is organized as follows. In \textcolor[rgb]{0.00,0.00,0.00}{Section 2}, we obtain global existence of weak solution for the system \eqref{eq: Landau-Lifshitz with Oseen-Frank} by using
the local well-posedness and blow-up results in the Appendix.
\textcolor[rgb]{0.00,0.00,0.00}{In Section 3}, we prove that the weak solution obtained in \textcolor[rgb]{0.00,0.00,0.00}{Section 2} is unique indeed.
At last, the local well-posedness and blow-up results for the Landau-Lifshitz system \eqref{eq: Landau-Lifshitz with Oseen-Frank} with general Oseen-Frank energy are obtained in \textcolor[rgb]{0.00,0.00,0.00}{the Appendix}.


\setcounter{equation}{0}
\section{Global existence of weak solutions in $\R^2$}

Let $E(t)=\int_{\R^2}W(\mn,\nabla\mn)(x,t)dx$ for $t\geq 0$ and $E_0=E(0)=\int_{\R^2}W(\mn_0,\nabla\mn_0)(x)dx$. Moreover,
\beno
E_R(\mn(\cdot,t);x)=\int_{B_R(x)}|\nabla\mn(y,t)|^2dy.
\eeno
For two constants $\tau$ and $T$ with $0\leq \tau<T$, we denote
$$
\begin{array}{ll}
V(\tau, T):&=\{\mn:\R^2\times [\tau,T]\rightarrow S^2|~\mn\mbox{ is measurable and satisfies }\\
&~~~~~~~\displaystyle{\rm{esssup}}_{\tau\le t\le T}\int_{\R^2}|\nabla \mn(\cdot,t)|^2dx+\int_{\tau}^{T}\int_{\R^2}|\nabla^2\mn|^2+|\partial_t
\mn|^2dxdt<\infty\}.
\end{array}
$$

\subsection{A priori estimates}
The following technical lemma \textcolor[rgb]{0.00,0.00,0.00}{can} be found in \cite{Struwe85}.
\begin{Lemma}\label{lem: struwe}
There are constants $C$ and $R_0$ such that for any $u\in V(0,T)$ and any $R\in(0,R_0]$, we have

\ben\label{eq:the embedding inequlity}
\int_{\R^2\times[0,T]}|\nabla u|^4dxdt&\le & C~{\rm esssup}_{0\le t\le T, x\in\R^2}\int_{B_{R}(x)}|\nabla u(\cdot,t)|^2dx\nonumber\\
&& \cdot\big(\int_{\R^2\times[0,T]}|\nabla^2 u|^2+R^{-2}\int_{\R^2\times[0,T]}|\nabla u|^2dxdt\big).
\een
\end{Lemma}

First of all,  we have the following basic energy estimates.
\begin{Lemma}[The basic energy estimates]\label{lem:energy estimate}
Assume that $\mn$ is a smooth solution of  the Landau-Lifshitz equation
 (\ref{eq: Landau-Lifshitz with Oseen-Frank}) in $(0,T)\times \R^2$ and the initial data $\mn_0\in H^1_b(\R^2)$.
Then, for all $0<t<T$ there holds
\beno
\int_{\R^2}W(\mn,\nabla\mn)(x,t)dx+\alpha\int_0^t\int_{\R^2}|\partial_t\mn|^2dxds\leq E_0.
\eeno
\end{Lemma}
{\it Proof:} Multiply $\partial_t\mn$ on both sides of the equation (\ref{eq: Landau-Lifshitz with Oseen-Frank}) and integrate on $\R^2$, then the
property $|\mn|=1$
implies that
\beno
\int_{\R^2}|\partial_t\mn|^2dx&=&\alpha \int_{\R^2} \partial_t\mn\cdot(\mn\times(\mh\times\mn))dx   +\beta\int_{\R^2}\partial_t\mn\cdot(\mn\times\mh)dx\\
&=& \alpha \int_{\R^2} \partial_t\mn\cdot\mh dx   +\beta\int_{\R^2}\partial_t\mn\cdot(\mn\times\mh)dx
\eeno
Noting that the definition of the molecular field $\mh$, we get
\beno
\frac{d}{d t}\int_{\R^2}W(\mn,\nabla\mn)(x,t)dx
= \int_{\R^2}(-\mh)\cdot \partial_t\mn dx.
\eeno
It follows that
\ben\label{eq:energy estimate1}
\int_{\R^2}|\partial_t\mn|^2dx+\alpha \frac{d}{d t}\int_{\R^2}W(\mn,\nabla\mn)(x,t)dx=\beta\int_{\R^2}\partial_t\mn\cdot(\mn\times\mh)dx.
\een
Now we estimate the term $\partial_t\mn\cdot(\mn\times\mh)$ as  in \cite{GuoHong1993}. The equation (\ref{eq: Landau-Lifshitz with
Oseen-Frank}) show that
\beno
\partial_t\mn=\alpha\mn\times (\mh\times \mn)+\beta\mn\times \mh,
\eeno
then we have
\beno
\mn\times\partial_t\mn=\alpha \mn\times\mh+\beta\mn\times(\mn\times\mh),
\eeno
hence using $\alpha^2+\beta^2=1$ we arrive at
\beno
\mn\times\partial_t\mn+\frac{\beta}{\alpha}\partial_t\mn=\frac{1}{\alpha}\mn\times\mh,
\eeno
which yields that
\ben\label{eq:partial t n n times h}
\partial_t\mn\cdot(\mn\times\mh)=\beta|\partial_t\mn|^2.
\een

Combining the estimates (\ref{eq:energy estimate1}) and (\ref{eq:partial t n n times h}), we have
\beno
\alpha^2\int_{\R^2}|\partial_t\mn|^2dx+\alpha \frac{d}{d t}\int_{\R^2}W(\mn,\nabla\mn)(x,t)dx=0,
\eeno
and the proof is completed by integrating with respect to time.
\endproof

As in \cite{Struwe85,GuoHong1993}, the key ingredient for global existence of weak solution is a local monotonicity inequality, and our results state as
follows.
\begin{Lemma}[The local monotonicity inequality]\label{lem:local monotonicity inequality}
Assume that $\mn$ is a smooth solution of  the Landau-Lifshitz equation (\ref{eq: Landau-Lifshitz with Oseen-Frank})
 in $(0,T)\times \R^2$ and the initial data $\mn_0\in H^1_b(\R^2)$.
Then, for all $0<t<T$ and $x_0\in \R^2$ there holds
\beno
E_R(\mn(\cdot,t);x_0)\leq E_{2R}(\mn_0(\cdot);x_0)+C_0\frac{t}{R^2}E_0,
\eeno
where $C_0$ is an absolute constant independent of $t,R$ and $\mn.$
\end{Lemma}
{\it Proof:} Let $\phi(x)$ be a smooth cut-off function satisfying $\phi(x)=1$ for $x\in B_R(x_0)$ and $\phi(x)=0$ when $|x-x_0|>2R.$
Multiply  $\partial_t\mn\phi^2$ on both sides of  (\ref{eq: Landau-Lifshitz with Oseen-Frank}), then we have
\beno
\int_{\R^2}|\partial_t\mn|^2\phi^2dx
= \alpha \int_{\R^2} \partial_t\mn\cdot\mh \phi^2dx   +\beta\int_{\R^2}\partial_t\mn\cdot(\mn\times\mh)\phi^2dx,
\eeno
and using the following relation
\beno
\frac{d}{d t}\int_{\R^2}W(\mn,\nabla\mn)(x,t)\phi^2(x)dx
= \int_{\R^2}(-\mh)\cdot \partial_t\mn \phi^2dx-2\int_{\R^2}W_{P_i^j}(\mn,\nabla\mn)\partial_t\mn^j\partial_i\phi \phi dx,
\eeno
hence we get
\beno
&&\int_{\R^2}|\partial_t\mn|^2\phi^2dx+\alpha\frac{d}{d t}\int_{\R^2}W(\mn,\nabla\mn)(x,t)\phi^2(x)dx\\
&\leq & 2\alpha\left|\int_{\R^2}W_{P_i^j}(\mn,\nabla\mn)\partial_t\mn^j\partial_i\phi \phi dx\right|+\beta\int_{\R^2}\partial_t\mn\cdot(\mn\times\mh)\phi^2dx,
\eeno
and using the equality of (\ref{eq:partial t n n times h}) for the term $\partial_t\mn\cdot(\mn\times\mh)$ again, there holds
 \beno
\frac{d}{d t}\int_{\R^2}W(\mn,\nabla\mn)(x,t)\phi^2(x)dx\leq  C(\alpha)\int_{\R^2}|\nabla\mn|^2|\nabla\phi|^2 dx\leq C_0\frac{1}{R^2}E_0.
 \eeno
Then the proof is complete.
\endproof

\begin{Lemma}[The positive diffusion]\label{lem:The positive diffusion}
Assume that $\mn$ is a smooth solution of  the Landau-Lifshitz equation  (\ref{eq: Landau-Lifshitz with Oseen-Frank})
in $(0,T)\times \R^2$ and the initial data $\mn_0\in H^1_b(\R^2)$.
Then there exists $\epsilon_{1}>0$, such that for all $R\in (0,R_0]$ with $R_0>0$, if
$$
{\rm esssup}_{0\le t\le T,x\in \R^2}\int_{B_{R}(x)}|\nabla\mn(\cdot,t)|^2dx<\epsilon_1,
$$
then there hold
\ben\label{eq:nabla2n estimate}
\int_{\R^2\times[0,T]}|\nabla^2\mn|^2dxdt\le C(1+TR^{-2})E_0,
\een
and
\ben\label{eq:nablan4 estimate}
\int_{\R^2\times[0,T]}|\nabla\mn|^4dxdt\le C\epsilon_{1}(1+TR^{-2})E_0.
\een
\end{Lemma}
{\it Proof:} Due to the embedding inequality (\ref{eq:the embedding inequlity}), it suffices to prove the first inequality (\ref{eq:nabla2n estimate}).
Since $$\frac{d}{dt}\int_{\R^2}W(\mn,\nabla\mn) dx=\int_{\R^2}\big(W_{n^{l}}-\nabla_{i}W_{p_{i}^{l}}\big)\cdot n_{t}^{l}dx=-\int_{\R^2}\mh\cdot\mn_t dx,$$
using the equation of (\ref{eq: Landau-Lifshitz with Oseen-Frank}) we have
\beno
\int_{\R^2}W(\mn,\nabla\mn)(\cdot,t) dx+\alpha\int_0^t\int_{\R^{2}}(\mn\times(\mh\times\mn))\cdot\mh dxds\leq E_0.
\eeno

Next we prove the positivity of the diffusion term.
Using Lemma \ref{lem:h decomposition} and $\mn\cdot\Delta\mn=-|\nabla\mn|^2$, we derive that
   \beno
   &&\int_{\R^2\times[0,T]}(\mn\times(\mh\times\mn))\cdot\mh dxdt\nonumber\\
   &&\ge \int_{R^2\times[0,T]}(\mn\times(\nabla_{i}W_{p_{i}^{l}}\times\mn))\nabla_{i}W_{p_{i}^{l}} dxdt
 -C\int_{\R^2\times[0,T]}|\nabla\mn|^2(|\nabla^2\mn|+|\nabla\mn|^2)dxdt\nonumber\\
&&\ge 2a\int_{\R^2\times[0,T]}\nabla_{i}W_{p_{i}^{l}} \cdot\Delta\mn dxdt
+2a\int_{\R^2\times[0,T]}\Delta\mn\cdot(\nabla_{i}W_{p_{i}^{l}}-2a\Delta\mn)dxdt\nonumber\\
   &&\quad+\int_{\R^2\times[0,T]}(\mn\times((\nabla_{i}W_{p_{i}^{l}}-2a\Delta\mn)\times\mn)) \cdot(\nabla_{i}W_{p_{i}^{l}}-2a\Delta\mn) dxdt\nonumber\\
    &&\quad-C\int_{\R^2\times[0,T]}|\nabla\mn|^2(|\nabla^2\mn|+|\nabla\mn|^2)dxdt\nonumber\\
   & &\ge  4a\int_{\R^2\times[0,T]}\big[a|\Delta\mn|^2+2(k_1-a)|\nabla{\rm div}\mn|^2+2(k_2-a)|\nabla(\nabla\times\mn\times\mn)|^2\nonumber\\
  &&\quad+2(k_3-a)|\nabla(\nabla\times\mn\cdot\mn)|^2\big] dxdt-C\int_{\R^2\times[0,T]}|\nabla\mn|^2(|\nabla^2\mn|+|\nabla\mn|^2)dxdt,\nonumber\\
&&\ge 3a^2\int_{\R^2\times[0,T]}|\Delta\mn|^2dxdt-C\int_{\R^2\times[0,T]}|\nabla\mn|^4dxdt,
   \eeno
and the first estimate (\ref{eq:nabla2n estimate}) follows from the embedding inequality (\ref{eq:the embedding inequlity}) by choosing a small
$\epsilon_1$.
\endproof

Concluding the above local monotonicity inequality in \textcolor[rgb]{0.00,0.00,0.00}{Lemma} \ref{lem:local monotonicity inequality} and the positive diffusion in \textcolor[rgb]{0.00,0.00,0.00}{Lemma} \ref{lem:The positive diffusion},
we have the following corollary.

\begin{Corollary} \label{cor:local bound of lapalace n}
Assume that $\mn$ is a smooth solution of  the Landau-Lifshitz equation (\ref{eq: Landau-Lifshitz with Oseen-Frank})
in $(0,T)\times \R^2$ and the initial data $\mn_0\in H^1_b(\R^2)$.
Then, there exists $R>0$ such that
$\sup_{x\in\R^2}E_{2R}(\mn_0(\cdot);x)\leq \frac{\epsilon_1}{2}$, and
\ben\label{eq:local bound of lapalace n}
\int_{\R^2\times[0,t]}|\nabla\mn|^4dxdt+
\int_{\R^2\times[0,t]}|\nabla^2\mn|^2dxdt\le C(E_0+\epsilon_1),
\een
hold for $t<\frac{\epsilon_1 R^2}{2C_0 E_0}$, where $C_0$ is given in Lemma \ref{lem:local monotonicity inequality}.
\end{Corollary}

Next, we use the idea of  Lemma \ref{lem:The positive diffusion} and the estimates in Corollary \ref{cor:local bound of lapalace n} to obtain a higher interior regularity of the solution.
\begin{Lemma} \label{lem:nabla 2 n}
Assume that $\mn$ is a smooth solution of  the Landau-Lifshitz equation (\ref{eq: Landau-Lifshitz with Oseen-Frank})
 in $(0,T)\times \R^2$ and the initial data $\mn_0\in H^1_b(\R^2)$.
Then  there is a constant $\epsilon_1$ such that for all $R\in(0,R_0]$, if
\beno
 {\rm esssup}_{0\le t\le T,x\in\R^2}\int_{B_{R}(x)}|\nabla\mn(\cdot,t)|^2dx<\epsilon_1,\eeno
then, for all $t\in(\tau, T)$ with $\tau\in (0,T)$, it holds that
\beno
\int_{\R^2}|\nabla^2\mn(\cdot,t)|^2dx+\int_{\tau}^t\int_{\R^2}|\nabla^3\mn(\cdot,s)|^2dxds \le C(\epsilon_1,E_0,\tau, T,\frac{T}{R^2}).
\eeno
\end{Lemma}
{\it Proof:}
First, we can differentiate $\nabla_\beta$ to (\ref{eq: Landau-Lifshitz with Oseen-Frank}),  multiply it by $\nabla_i \mh$($i=1,2$), and we get
\begin{align}\label{eq:nabla 2 n}
&\frac{d}{dt}\int_{\R^2}a|\Delta \mn|^2+(k_1-a)|\nabla{\rm div}\mn|^2
+(k_2-a)|\nabla(\mn\times(\nabla\times\mn))|^2dx\nonumber\\
&+\frac{d}{dt}\int_{\R^2}(k_3-a)|\nabla(\mn\cdot(\nabla\times\mn))|^2dx\nonumber\\
\leq &-\alpha\int_{\R^2}\nabla_{i}(\mn\times(\mh\times\mn))\cdot\nabla_{i}\mh dx
-\beta \int_{\R^2}\nabla_{i}(\mn\times\mh)\cdot\nabla_{i}\mh dx\nonumber\\
&+C\int_{\R^2}[|\nabla\mn_t| |\nabla\mn||\nabla^2\mn|+|\mn_t| |\nabla^2\mn|^2]dx\nonumber\\
\leq& -\alpha\int_{\R^2}(\mn\times(\nabla_{i}\mh\times\mn))\cdot\nabla_{i}\mh dx\nonumber\\
&+C\int_{\R^2}[|\nabla\mn_t| |\nabla\mn||\nabla^2\mn|+|\mn_t| |\nabla^2\mn||\nabla\mn|^2+|\nabla^2\mn|^2|\nabla\mn|^2+|\nabla\mn||\nabla^3\mn||\nabla^2\mn|]dx\nonumber\\
\leq& -\alpha\int_{\R^2}(\mn\times(\nabla_{i}\mh\times\mn))\cdot\nabla_{i}\mh dx+C\int_{\R^2}[|\nabla^2\mn|^2|\nabla\mn|^2+|\nabla\mn||\nabla^3\mn||\nabla^2\mn|]dx
\end{align}

Note the fact that $|\mn\cdot\nabla_{i}\Delta\mn|\leq C|\nabla\mn||\nabla^2\mn|$, and similar estimates as in Lemma \ref{lem:The positive diffusion} imply
\beno
\int_{\R^2}(\mn\times(\nabla_{i}\mh\times\mn))\cdot\nabla_{i}\mh dx\geq 3 a^2\int_{\R^2}|\nabla^3\mn|^2-C\int_{\R^2}|\nabla\mn|^2|\nabla^2\mn|^2dx.
\eeno

Due to the interpolation inequality
\beno
\|\nabla^2\mn\|_4\leq C \|\nabla^2\mn\|_2^{1/2} \|\nabla^3\mn\|_2^{1/2},
\eeno
we have
\beno
\int_{\R^2}|\nabla\mn|^2|\nabla^2\mn|^2dx\leq \delta \|\nabla^3\mn\|_2^2+C(\delta)\|\nabla\mn\|_4^4\|\nabla^2\mn\|_2^2,
\eeno
thus Gronwall's inequality and Corollary \ref{cor:local bound of lapalace n} imply the required estimates.
\endproof


Indeed, using the above idea by induction, one can prove the smooth property of $\mn$, and we omit the proof (similar arguments for  Ericksen-Leslie system, see \cite[Corollary 4.6]{WangW2014}).
\begin{Corollary} \label{cor:higher regularity}
Assume that $\mn$ is a smooth solution of  the Landau-Lifshitz equation (\ref{eq: Landau-Lifshitz with Oseen-Frank})  in $(0,T)\times \R^2$ and the initial data $\mn_0\in H^1_b(\R^2)$. Then there is a constant $\epsilon_1>0$ such that for all $R\in(0,R_0]$, if
\beno
 {\rm esssup}_{0\le t\le T,x\in\R^2}\int_{B_{R}(x)}|\nabla\mn(\cdot,t)|^2dx<\epsilon_1,\eeno
then, for all $t\in(\tau, T)$ with $\tau\in (0,T)$, for any $l\geq 1$ it holds that
\ben\label{eq:l regularity}
&&\int_{\R^2}|\nabla^{l+1}\mn(\cdot,t)|^2dx +\int_{\tau}^t\int_{\R^2}|\nabla^{l+2}\mn(\cdot,s)|^2dxds
\nonumber\\
&&\le C(l,\epsilon_1,E_0,\tau, T, \frac{T}{R^2}).
\een
Moreover, $\mn$ is regular for all $t\in(0,T)$.
\end{Corollary}

\subsection{Existence of global weak solution}

Now we complete the proof of the existence part in Theorem \ref{thm:global}. Similar to \cite{Struwe85,LLW,WangW2014}, we sketch its step for completeness.

For any data $\mn_0\in H_{b}^{1}(\R^2;S^2)$, one can approximate it by a sequence of smooth maps $\mn_0^k$ in $H_{b}^{1}(\R^2;S^2)$, and we can assume
that $\nabla\mn_0^k\in H_{b}^{4}(\R^2;S^2)$ (see \cite{SchoenU82}).
Due to the absolute continuity property of the integral, for any $\epsilon_1>0$, there exists $R_0\ge R_1>0$ such that
$$
 \sup_{x\in \R^2}\int_{B_{R_1}(x)}|\nabla\mn_0|^2dx\le \epsilon_1,
$$
and by the strong convergence of $\mn_0^k$,
$$
 \sup_{x\in \R^2}\int_{B_{R_1}(x)}|\nabla\mn_0^k|^2dx\le 2\epsilon_1
$$
for a sufficient large $k$. Without loss of generality, we assume that it holds for all $k\geq 1.$

For the data $\mn_0^k$, by Theorem \ref{thm:local existence} there exists a time $T^k$ and a strong solution $\mn^k$  such that
$$
\nabla\mn^k\in C\left([0,T^k];H^{4}(\R^2)\right).
$$
Hence there exists $T_0^k\leq T^k$ such that
$$
 \sup_{0<t<T_0^k,\\ x\in \R^2}\int_{B_{R}(x)}|\nabla\mn^k(y,t)|^2dy\leq (8+\frac1a)\epsilon_1,
$$
where $R\leq R_1/2.$
However, by the local monotonic inequality in Lemma \ref{lem:local monotonicity inequality}, we have $T_0^k\geq
\frac{\epsilon_1R_1^2}{4C_0E_0}=T_0>0$  uniformly.
For any $0<\tau<T_0$, by the estimates in Corollary \ref{cor:higher regularity} for any $l\geq 1$ we get
\begin{eqnarray}\label{eq:5.1}
\sup_{\tau<t<T_0}\int_{\R^2} |\nabla^{l+1}\mn^k|^2(\cdot,t)dx +\int_{\tau}^{T_0}\int_{\R^2} |\nabla^{l+2}\mn^k(\cdot,s)|^2dxds \le C(l,\epsilon_1,E_0,\tau,  T_0, \frac{T_0}{R^2}).
\end{eqnarray}
Moreover, the energy inequality in Lemma \ref{lem:energy estimate}, {\it a priori}  estimates in Lemma \ref{lem:The positive diffusion} and the equation
(\ref{eq: Landau-Lifshitz with Oseen-Frank}) yield that
\ben\label{eq:5.2}
E(\mn^k)(t)\leq E_0,\quad  0<t<T^k,
\een
and
\begin{eqnarray}\label{eq:5.3}
\int_{\R^2\times [0,T_0^k]}\big(|\nabla^2 \mn^k|^2+|\partial_t \mn^k|^2+|\nabla \mn^k|^4\big)dxdt\leq C(\epsilon_1,C_0, E_0).
\end{eqnarray}

Hence the above estimates (\ref{eq:5.1})-(\ref{eq:5.3}) and Aubin-Lions Lemma yield that there exists a solution $\mn-b\in W^{2,1}_2 (\R^2\times
[0,T_0];\R^3)$ such that
(at most up to a subsequence)
\begin{eqnarray*}
&&\mn^k-b\rightarrow \mn-b,\quad {\rm locally\,\, in}\quad W^{3,1}_2(\R^2\times (0,T_0);\R^3).
\end{eqnarray*}
By (\ref{eq:5.2}), $\nabla\mn(t)\rightharpoonup \nabla\mn_0$ weakly in $L^2(\R^2)$, thus
$E(\mn_0)\leq\liminf_{t\rightarrow0} E(\mn(t)).$ On the other hand, by the energy estimates of $(\mn^k)$, we have
$$E(\mn_0)\geq\limsup_{t\rightarrow0} E(\mn(t)).$$
Hence,  $\nabla\mn(t)\rightarrow \nabla\mn_0$ strongly in $L^2(\R^2)$ and $\mn$ is the solution of the equation (\ref{eq: Landau-Lifshitz with
Oseen-Frank}) with the initial data $\mn_0.$ From the weak limit of regular estimates (\ref{eq:5.1}), we know that $\mn\in C^{\infty}(\R^2\times(0,T_0])$ and $\nabla^{l+1}\mn(\cdot ,T_0)\in L^2(\R^2) $ for any $l\geq1$. By Theorem \ref{thm:local existence}, there exists a unique smooth
solution of (\ref{eq: Landau-Lifshitz with Oseen-Frank}) with the initial data $\mn(\cdot,T_0)$, which is still written as $\mn$, and blow-up criterion
yields that if $\mn$ blows up at finite time $T^*$, then
\beno
\|\nabla\mn\|_{L^{\infty}(\R^2)}(t)\rightarrow\infty,\quad {\rm as} \quad t\rightarrow T^*.
\eeno
As a result, we have
\ben\label{eq:blow-up}
|\nabla^{4}\mn|(x,t)\not\in L^{\infty}_tL^{2}_x((T_0,T^*)\times\R^2)
\een
We assume that $T_1$ is the first singular time of  $\mn$, then we have
\begin{eqnarray*}
\mn\in C^{\infty}(\R^2\times (0,T_1);  S^2)\quad {\rm and}
\quad \mn\not\in C^{\infty}(\R^2\times (0,T_1]; S^2);
\end{eqnarray*}
and by Corollary \ref{cor:higher regularity} and (\ref{eq:blow-up}), there exists $\epsilon_0>0$ such that
\begin{eqnarray*}
\lim\sup_{t\uparrow T_1}\sup_{x\in \R^2}\int_{B_R(x)}|\nabla\mn|^2(\cdot,t)\geq \epsilon_0,\quad \forall R>0.
\end{eqnarray*}
Finally, since $\mn-b\in C^0([0,T_1], L^2(\R^2))$ by the interpolation inequality (similarly see   P330, \cite{LLW}), we can define
$$ \mn(T_1)-b=\lim_{t\uparrow T_1} \mn(t)-b\quad {\rm in}\quad  L^2(\R^2).$$
On the other hand, by the energy inequality $\nabla\mn\in L^{\infty}(0,T_1;L^2(R^2))$, hence
$\nabla\mn(t)\rightharpoonup \nabla\mn(T_1)$. Similarly we can extend $T_1$ to $T_2$ and so on. It's easy to check that the energy loss at every singular
time $T_i$ for $i\geq 1$ is at least $\epsilon_1$, thus the number $L$ of the singular time is finite. Moreover, singular points at every singular time are finite by similar arguments as in \cite{Struwe85}, since $\partial_tu\in L^2_{x,t}$ in Lemma \ref{lem:energy estimate} and the local monotonicity inequality in Lemma \ref{lem:local monotonicity inequality} hold. Assume that singular points are $(x_i^{j},T_i)$ with $1\leq j\leq L_i$ and $i\leq L$, and we have
\begin{eqnarray*}
\lim\sup_{t\uparrow T_i}\int_{B_R(x_i^{j})}|\nabla\mn|^2(\cdot,t)\geq \epsilon_0,\quad \forall R>0.
\end{eqnarray*}
The proof is complete.\endproof


\setcounter{equation}{0}
\section{Uniqueness of weak solution}

In this section, we follow the same route as in \cite{WangWZ2013} and prove the following uniqueness theorem. The main difference is to deal with the Schr\"{o}dinger part $\mn\times\mh$.
\begin{Theorem}\label{thm:uniqueness}
 Let $\mn^1$ and $\mn^2$ be two weak solutions of the Landau-Lifshitz equation (\ref{eq: Landau-Lifshitz with Oseen-Frank}) in $\R^2$ obtained in Theorem
 \ref{thm:global}
 with the same initial data $\mn_0$. Then we have
 $$\mn^1(t)=\mn^2(t)$$
 for any $t\in [0,+\infty)$.
\end{Theorem}

Let $\mn^1$ and $\mn^2$ be two weak solutions of  the Landau-Lifshitz equation (\ref{eq: Landau-Lifshitz with Oseen-Frank}) in $\R^2$ obtained in Theorem
\ref{thm:global}
 with the same initial data $\mn_0$. Let
\beno
\delta_{\mn}=\mn^2-\mn^1,
\eeno
then we infer that
\ben\label{eq:delta n}
\partial_t\delta_{\mn}=\alpha\delta_{\mn\times(\mh\times\mn)}+\beta\delta_{\mn\times\mh}.
\een
Here and in what follows, we denote $f^i=f(\mn^i)$ for $i=1,2$ and $\delta_f=f^2-f^1$ if $f$ is a function of $\mn$.



Different with \cite{Struwe85,GuoHong1993}, it's not easy to obtain the positivity of the diffusion term $\nabla\delta\mn$ under the metric of $L^2$ norm, since we can't use the property of $\triangle\mn\cdot\mn=-|\nabla\mn|^2$ from $|\mn|=1.$ Instead, we
 introduce a type of weak Oseen-Frank  metric
\beno
W(t)=\sup_{j\ge 0}2^{-2js}\int_{\R^2}W^j(t,x)dx+\|\Delta_{-1}\delta_{\mn}\|_{2}^2
\eeno
with $s\in (0,1)$ and
\begin{align}
W^j(x,t)=&a|\nabla\Delta_j\delta_{\mn}|^2+(k_1-a)|{\rm div}\Delta_j\delta_{\mn}|^2\nonumber\\
&+(k_2-a)|\mn^2\times({\mathbf{\nabla}}\times\Delta_j\delta_{\mn})|^2+(k_3-a)|\mn^2\cdot(\mathbf{\nabla}\times\Delta_j\delta_{\mn})|^2.\nonumber
\end{align}

The proof of Theorem \ref{thm:uniqueness} is based on the following two propositions. To state them, we introduce
\beno
\bar{h}(t)\eqdef 1+\|(\nabla\mn^1,\nabla\mn^2)\|_{4}^4+\|(\partial_t\mn^1,\partial_t\mn^2)\|_{2}^2+\|(\nabla\mn^1,\nabla\mn^2)\|_{H^1}^2.
\eeno

\begin{Proposition}\label{prop:n-low}
It holds that
\begin{align*}
\frac{d}{dt}\|\Delta_{-1}\delta_{\mn}\|_{2}^2
\le C\bar{h}(t)W(t).
\end{align*}
\end{Proposition}

\begin{Proposition}\label{prop:W}
For  any $j\geq 0$ and $\epsilon>0$, it holds that
\begin{align*}
\frac{d}{dt}\int_{\R^2}W^j(x,t)dx+3\alpha a^2\|\Delta_j\na^2\delta_\mn\|_2^2\leq
C2^{2js}\bar{h}(t)W(t)
+\epsilon \sum_{l=j-9}^{j+9}2^{4l}\|\Delta_{l}\delta_{\mn}\|_{2}^2.
\end{align*}
\end{Proposition}

For the moment, let us assume that these propositions are correct and complete the proof of Theorem \ref{thm:uniqueness}.
Assume that $T_1^i$ is the first blow-up time of $\mn^i$ with $i=1,2$. We know from Lemma \ref{lem:The positive diffusion} that
\ben\label{eq:solu-R1}
\int_{\R^2\times[0,T_1-\theta]}|\nabla^2\mn^i|^2+|\nabla\mn^i|^4dxdt<+\infty,
\een
where $\theta>0$ and $T_1=\min\{T^1_1,T^2_1\}$. And using the equation (\ref{eq: Landau-Lifshitz with Oseen-Frank}), we get
\ben\label{eq:solu-R2}
\partial_t\mn^i\in L^2((0,T_1-\theta)\times\R^2).
\een
Proposition \ref{prop:n-low} and Proposition \ref{prop:W} ensure that
\begin{align*}
\f d {dt}\Big(\int_{\R^2}W^jdx+\|\Delta_{-1}\delta_\mn\|_2^2\Big)
+ca2^{4j}\|\Delta_j\delta_\mn\|_2^2
\le C2^{2js}\bar{h}(t)W(t)
+\epsilon \sum_{l=j-9}^{j+9}2^{4l}\|\Delta_{l}\delta_{\mn}\|_{2}^2.
\end{align*}
Noting that $\int_{\R^2}W^jdx+\|\Delta_{-1}\delta_\mn\|_2^2\ge c2^j\|\Delta_j\delta_\mn\|_2^2$, we deduce by taking $\epsilon$ small enough that
\beno
W(t)\le C\int_0^t\bar{h}(\tau)W(\tau)d\tau.
\eeno
By (\ref{eq:solu-R1}) and (\ref{eq:solu-R2}), $\bar{h}(t)\in L^1(0,T_1-\theta)$. Then by Gronwall's inequality, we get $W(t)=0$ for $t\in [0,T_1-\theta]$
for any $\theta>0$.
Hence, $\mn^1(t)=\mn^2(t)$ on $[0,T_1)$ with $T_1>0$ the first singular time of the solution $\mn^1$ or $\mn^2$.
Since $\mn^i\in C_w([0,+\infty); H^1_b)$, $\mn^1(T_1)=\mn^2(T_1)$.
Then the same arguments show that there exists a $T_2>T_1$ such that $\mn^1(t)=\mn^2(t)$ on $[T_1,T_2)$,
where $T_2$ is the second singular time of the solution $\mn^1$ or $\mn^2$.
Since the number of singular time is finite,
we can conclude that $\mn^1(t)=\mn^2(t)$ for $t\in [0,+\infty)$. \endproof

\subsection{Littlewood-Paley theory and nonlinear estimates}

 Let us recall some basic facts on Littlewood-Paley theory (see \cite{Che}  for more details). Choose two nonnegative radial functions $\chi,\phi\in {\cal
 S}(\R^n)$ supported respectively in
$\{\xi\in \R^n,|\xi|\le \frac{4}{3}\}$ and $\{\xi\in \R^n, \frac{3}{4}\le |\xi|\le \frac{8}{3}\}$ such that for any $\xi\in \R^n$,
$$
\chi(\xi)+\sum_{j\ge 0}\phi(2^{-j}\xi)=1.
$$
The frequency localization operator $\Delta_j$ and $S_j$ are defined by
\beno
&&\Delta_j f=\phi(2^{-j}D)f=2^{nj}\int_{\R^n}h(2^{j}y)f(x-y)dy,~~~~\mbox{for}~~ j\ge 0,\\
&&S_j f=\chi(2^{-j}D)f=\sum_{-1\le k\le j-1}\Delta_{k}f=2^{nj}\int_{\R^n}\tilde{h}(2^{j}y)f(x-y)dy,\\
&&\Delta_{-1}f=S_0f, ~~\Delta_j f=0~~ \mbox{for}~~ j\le -2,
\eeno
where $h={\cal F}^{-1}\phi$, $\tilde{h}={\cal F}^{-1}\chi$.  With this choice of $\phi$, it is easy to verify that
\ben\label{eq:A.1}
\Delta_j\Delta_kf=0,~~ {\rm if}~~|j-k|\geq 2;~~\Delta_j(S_{k-1}f\Delta_kf)=0,~~ {\rm if}~~|j-k|\geq 5.
\een

In terms of $\Delta_j$, the norm of the inhomogeneous Besov space $B^s_{p,q}$ for $s\in \R,$ and $p,q\geq 1$ is defined by
\beno
\|f\|_{B^s_{p,q}}\eqdef\big\|\{2^{js}\|\Delta_jf\|_p\}_{j\ge -1}\big\|_{\ell^q},
\eeno
and
\beno
\|f\|_{B^s_{p,\infty}}\eqdef\sup_{j\ge -1}\{2^{js}\|\Delta_jf\|_p\}.
\eeno


We will constantly use the following Bernstein's inequality \cite{Che}.
\begin{Lemma}\label{lem:Berstein}
Let $c\in (0,1)$ and $R>0$. Assume that $1\leq p\leq q\leq \infty$ and $f\in L^p(\R^n)$. Then
\beno
&&{\rm supp} \hat{f}\subset\big\{|\xi|\leq R\big\}\Rightarrow \|\partial^{\alpha}f\|_{q}\leq CR^{|\alpha|+n(\f1p-\f1q)}\|f\|_{p},\label{eq:A.3}\\
&&{\rm supp} \hat{f}\subset\big\{cR \leq |\xi|\leq R\big\}\Rightarrow \|f\|_{p}\leq
CR^{-|\alpha|}\sup_{|\beta|=|\alpha|}\|\partial^{\beta}f\|_{p},\label{eq:A.4}
\eeno
where the constant $C$ is independent of $f$ and $R$.
\end{Lemma}

We need the following nonlinear estimates, seeing \cite{WangWZ2013} for more details.

\begin{Lemma}\label{lem:product3-1}
Let $s\in (0,1)$. For any $j\ge -1$, we have
\beno
\|\Delta_j(fgh)\|_2\le C2^{js}\big(\|f\|_\infty+\|\na f\|_2\big)\|g\|_{B^{1-s}_{2,\infty}}\|h\|_2.
\eeno
\end{Lemma}

\begin{Lemma}\label{lem:product3-2}
Let $s\in (0,1)$. For any $j\ge -1$, we have
\begin{align*}
\|\Delta_j(f\na gh)\|_2\le& C2^{js}\|g\|_{B^{1-s}_{2,\infty}}\big(\|f\|_\infty\|h\|_{H^1}+\|\na f\|_{4}\|h\|_{4}+\|\na^2f\|_2\|h\|_2\big)\\
&+C2^{\frac{js}{2}}\|f\|_\infty\|h\|_{4}\|g\|^{\f12}_{B^{1-s}_{2,\infty}}\sum_{l=j-9}^{j+9}2^{\frac{l}{2}}
\|\Delta_l\nabla g\|_{2}^{\frac{1}{2}}.
\end{align*}
\end{Lemma}

\begin{Lemma}\label{commutator}
Let $s\in (0,1)$. For any $j\ge -1$, it holds that
\beno
\big\|[\Delta_j,f]\na g\big\|_{2}\le C 2^{\frac{js}{2}}\|\nabla f\|_{4}\|g\|_{B^{-s}_{2,\infty}}^\f12\sum_{|j'-j|\le 4}2^{\frac{j'}{2}}\|\Delta_{j'}g\|_{2}^{\frac{1}{2}}
+C2^{js}\|g\|_{B^{-s}_{2,\infty}}\big(\|f\|_\infty+\|\nabla^2 f\|_{2}\big).
\eeno
\end{Lemma}

\subsection{Proof of Proposition \ref{prop:n-low}}

Using the equation (\ref{eq:delta n}), we obtain
\beno
\f12\frac{d}{dt}\|\Delta_{-1}\delta_{\mn}\|_{2}^2=
\alpha\big\langle\Delta_{-1} \delta_{\mn\times(\mh\times\mn)}, \Delta_{-1}\delta_{\mn}\big\rangle +\beta\big\langle\Delta_{-1}\delta_{\mn\times\mh},
\Delta_{-1}\delta_{\mn}\big\rangle
\triangleq I.
\eeno
Recall that the formula of $\mh$ in (\ref{p2}),
and we could write $I$ as
\begin{align*}
I=&\big\langle\Delta_{-1}(\cM\nabla^{2}\delta_{\mn}),\Delta_{-1}\delta_{\mn}\big\rangle
+\sum_{i=1,2}\big\langle\Delta_{-1}(\cM\nabla^2\mn^i\delta_{\mn}),\Delta_{-1}\delta_{\mn}\big\rangle\nonumber\\
&+\sum_{i=1,2}\big\langle\Delta_{-1}(\cM\nabla\mn^i\delta_{\nabla\mn}),\Delta_{-1}\delta_{\mn}\big\rangle
+\sum_{i,k=1,2}\big\langle\Delta_{-1}(\cM\nabla\mn^i\nabla\mn^k\delta_{\mn}),\Delta_{-1}\delta_{\mn}\big\rangle.\\
=& I_1+\cdots+I_4.
\end{align*}
Here and in what follows, we denote by $\cM$ a polynomial function of $(\mn^1,\mn^2)$ with degree no greater than  $4$,
which may be different from line to line.
Then by Lemma \ref{lem:product3-1} (for $I_2,I_4$), Lemma \ref{lem:product3-2} and Lemma \ref{lem:Berstein} (for $I_1,I_3$), we get
\beno
|I|\le C\bar{h}(t)W(t).
\eeno
Thus the proof is complete. \endproof

\subsection{Proof of Proposition \ref{prop:W}}

Let us first derive the following evolution inequality for the Oseen-Frank density.
\begin{Lemma}\label{lem:W-evolution}
For any $j\geq 0$, it holds that
\begin{align*}
\frac{d}{dt}\int_{\R^2}W^j(t,x)dx+3\alpha a^2\|\Delta_j\na^2\delta_\mn\|_2^2\leq B_1+\cdots+ B_6,
\end{align*}
where $B_i$ will be given in the proof.
\end{Lemma}

The key part of the above lemma is the positivity of the diffusion term $ \mn\times(\triangle_j\delta_{\mh}\times\mn))\cdot\triangle_j\delta_{\mh}$. It's important to analysis the main parts of $\triangle_j\delta_{\mh}$ and $\triangle_j\delta_{\mn\times(\mh\times\mn)}$ (the second derivative terms).
Using $\mcurl(f u)=f\mcurl u+\nabla f\times u$, $\mh$ in (\ref{p2}) can be rewritten as
\ben\label{eq:h vertical n}
\mh&=&2a\Delta\mn+2(k_1-a)\nabla{\rm div}\mn-2(k_2-a)\mcurl\mcurl\mn-2(k_3-k_2)(\nabla\mcurl\mn\cdot\mn)\times\mn \nonumber\\
&&-2(k_3-k_2)\big(2(\mn\cdot\mcurl\mn)\mcurl\mn+(\nabla\mn\cdot\mcurl\mn)\times\mn\big),
\een
hence the main parts of $\Delta_{j}\delta_{\mh}$ is
\ben\label{eq:W1}
W_1&=&2a\Delta\Delta_j\delta_{\mn}+2(k_1-a)\nabla{\rm div}\Delta_j\delta_{\mn}-2(k_2-a)\mcurl\mcurl\Delta_{j}\delta_{\mn}\nonumber\\
&&-2(k_3-k_2)(\nabla\mcurl\Delta_j\delta_{\mn}\cdot\mn^2)\times\mn^2.
\een
Note that by (\ref{eq:h vertical n})
\begin{align*}
&(\mh\cdot\mn)\mn\nonumber\\
&=\left(-2a|\nabla\mn|^2+2(k_1-a)\mn\cdot\nabla{\rm div}\mn-2(k_2-a)\mn\cdot\mcurl\mcurl\mn-4(k_3-k_2)(\mn\cdot\mcurl\mn)^2\right)\mn\nonumber\\
&=2(k_1-a)(\mn\cdot\nabla{\rm div}\mn)\mn-2(k_2-a)(\mn\cdot\mcurl\mcurl\mn)\mn\nonumber\\
&\quad-2a|\nabla\mn|^2\mn-4(k_3-k_2)(\mn\cdot\mcurl\mn)^2\mn,
\end{align*}
and $\mn\times(\mh\times\mn)=\mh-(\mh\cdot\mn)\mn$. We deduce
\ben\label{eq:hn-d}
&&\delta_{\mn\times(\mh\times\mn)}\nonumber\\
&&=2a\Delta\delta_{\mn}+2(k_1-a)\nabla{\rm
div}\delta_{\mn}-2(k_2-a)\mcurl\mcurl\delta_{\mn}-2(k_3-k_2)(\nabla\mcurl\delta_{\mn}\cdot\mn^2)\times\mn^2\nonumber\\
&&\quad-2(k_1-a)(\mn^2\cdot\nabla{\rm div}\delta_{\mn})\mn^2+2(k_2-a)(\mn^2\cdot\mcurl\mcurl\delta_{\mn})\mn^2\nonumber\\
&&\quad+\sum_{i=1,2}(\cM\delta_{\mn}\nabla^2\mn^i+\cM\delta_{\nabla\mn}\nabla\mn^i)+\sum_{i,k=1,2}\cM\nabla\mn^i\nabla\mn^k\delta_{\mn}.
\een
Denote the main parts of  $\triangle_j\delta_{\mn\times(\mh\times\mn)}$ as follows.
\ben\label{eq:H1}
H_1&=&2a\Delta\Delta_j\delta_{\mn}+2(k_1-a)\nabla{\rm
div}\Delta_j\delta_{\mn}-2(k_2-a)\mcurl\mcurl\Delta_j\delta_{\mn}\nonumber\\
&&-2(k_3-k_2)(\nabla\mcurl\Delta_j\delta_{\mn}\cdot\mn^2)\times\mn^2
-2(k_1-a)(\mn^2\cdot\nabla{\rm div}\Delta_j\delta_{\mn})\mn^2\nonumber\\
&&+2(k_2-a)(\mn^2\cdot\mcurl\mcurl\Delta_j\delta_{\mn})\mn^2.
\een
\begin{Lemma}\label{lem:W1 H1 positivity}
Assume that $W_1,H_1$ state as above, then we have
\begin{align*}
\f14\int_{\R^2}W_1\cdot H_1dx
\ge\frac{3}{4}a^2\|\Delta\Delta_j\delta_{\mn}\|_{2}^2-B_1.
\end{align*}
where
\beno
B_1=|\langle\cM\nabla\mn^2\Delta_{j}\nabla^2\delta_{\mn},\Delta_{j}\nabla\delta_{\mn}\rangle|.
\eeno
\end{Lemma}
{\bf Proof of Lemma \ref{lem:W1 H1 positivity}.}
Let $S_1=\Delta_j\Delta\delta_\mn$, and
\beno
 H_2=(k_1-a)\nabla{\rm
div}\Delta_j\delta_{\mn}-(k_2-a)\mcurl\mcurl\Delta_j\delta_{\mn}-(k_3-k_2)(\nabla\mcurl\Delta_j\delta_{\mn}\cdot\mn^2)\times\mn^2.
\eeno
Then we find
\begin{align*}
\f14\int_{\R^2}W_1\cdot H_1dx=&\int_{\R^2}\big(aS_{1}+H_2\big)\cdot\big(aS_{1}+H_2-(\mn^2\cdot H_2)\mn^2\big)dx\nonumber\\
=&a^2\|S_1\|_{2}^2+a\big\langle H_2, S_1\big\rangle+\|H_2\times\mn^2\|_2^2+a\big\langle S_1, \mn^2\times(H_2\times\mn^2)\big\rangle\nonumber\\
&+\frac{a^2}{4}\|\mn^2\times S_1\|_2^2-\frac{a^2}{4}\|\mn^2\times S_1\|_2^2\nonumber\\
\ge& \f34 a^2\|\Delta_{j}\Delta\delta_{\mn}\|_{2}^{2}+a\big\langle H_2,S_1\big\rangle.
\end{align*}
Furthermore, by Lemma \ref{lem:positive of h} we have
\begin{align*}
\big\langle H_2, S_1\big\rangle=&(k_1-a)\|\nabla{\rm div}\Delta_{j}\delta_{\mn}\|_{2}^{2}+(k_2-a)\|\nabla\mcurl\Delta_{j}\delta_{\mn}\|_{2}^{2}\nonumber\\
&+(k_3-k_2)\big\langle\nabla\Delta_{j}\mcurl\delta_{\mn}\cdot\mn^2,\nabla\Delta_{j}\mcurl\delta_{\mn}\cdot\mn^2\big\rangle-B_1\geq -B_1.
\end{align*}
Hence we get
\begin{align*}
\f14\int_{\R^2}W_1\cdot H_1dx
\ge\frac{3}{4}a^2\|\Delta\Delta_j\delta_{\mn}\|_{2}^2-B_1.
\end{align*}
The proof is complete. \endproof

\no{\bf Proof of Lemma \ref{lem:W-evolution}.}\, Due to the definition of $W^j$, we have
\begin{align*}
\frac{d}{dt}\int_{\R^2}{W}^{j}(t,x)dx
&=\int_{\R^2}-\nabla_i{W}_{p_{i}^{k}}^{j}\partial_{t}\Delta_{j}\delta_{n_k}+{W}_{n_l}^{j}(n_{l}^{2})_{t}dx\nonumber\\
&\triangleq -\int_{\R^2}\nabla_i{W}_{p_{i}^{k}}^{j}\partial_{t}\Delta_{j}\delta_{n_k}dx+B_1.
\end{align*}
Using the equation (\ref{eq: Landau-Lifshitz with Oseen-Frank}), we get
\begin{align*}
-\int_{\R^2}(\nabla_{i}{W}^j_{p_{i}^{l}})\Delta_j\delta_{\mn_{t}}dx
=&-\alpha\int_{\R^2}(\nabla_{i}{W}^j_{p_{i}^{l}})\Delta_j\big(\mn^2\times(\mh^2\times\mn^2)
-\mn^1\times(\mh^1\times\mn^1)\big)dx\nonumber\\
&-\beta\int_{\R^2}(\nabla_{i}{W}^j_{p_{i}^{l}})\Delta_j\big(\mn^2\times\mh^2
-\mn^1\times\mh^1\big)dx\nonumber\\
\triangleq& \alpha I'+\beta I''.
\end{align*}
So, we conclude that
\ben\label{eq:bar w}
&&\frac{d}{dt}\int_{\R^2}{W}^{j}(t,x)dx\leq \alpha I'+\beta I''+B_2.
\een

As in Lemma \ref{lem:h decomposition}, we have
\begin{align*}
\nabla_{\alpha}{W}_{p_{\alpha}^{l}}^j=&2a\Delta_j\Delta\delta_{\mn}+2(k_1-a)\nabla{\rm
div}\Delta_{j}\delta_{\mn}-2(k_2-a)\mcurl\mcurl\Delta_{j}\delta_{\mn}\nonumber\\
&-2(k_3-k_2)\mcurl((\mn^2\cdot\mcurl\Delta_{j}\delta_{\mn})\mn^2),
\end{align*}
where $p_{\alpha}^{l}=\nabla_{\alpha}\Delta_j(\mn^2-\mn^1)_l$, and we get
\ben\label{eq:nabla W j}
\nabla_{\alpha}{W}^j_{p_{\alpha}^{l}}
=W_1+\cM\nabla\mn^2\nabla\Delta_j\delta_{\mn},
\een
that is, $W_1$ is also the main part of $\nabla_{\alpha}{W}^j_{p_{\alpha}^{l}}.$
Then we have
\begin{align*}
I''=&-\int_{\R^2}(\nabla_{i}{W}^j_{p_{i}^{l}})\cdot\big(\mn^2\times\Delta_j\delta_{\mh}\big)
dx-\int_{\R^2}(\nabla_{i}{W}^j_{p_{i}^{l}})\cdot\big([\Delta_j,\mn^2\times]\delta_{\mh}\big)
dx\\
&-\int_{\R^2}(\nabla_{i}{W}^j_{p_{i}^{l}})\Delta_j\big(\delta_{\mn}\times\mh^1\big)
dx\\
\triangleq& B_5+B_3+B_4,
\end{align*}
where $B_5$ can be further decomposed into
\begin{align*}
B_5=&-\int_{\R^2}(\nabla_{i}{W}^j_{p_{i}^{l}})\cdot\big(\mn^2\times(\Delta_j\delta_{\mh}-W_1)\big)dx
-\int_{\R^2}(\nabla_{i}{W}^j_{p_{i}^{l}}-W_1)\cdot\big(\mn^2\times W_1\big).
\end{align*}

On the other hand, for the estimate of $I',$ we have
\begin{align*}
I'=&-\int_{\R^2}W_1\cdot H_1dx-\int_{\R^2}\big(\nabla_{i}{W}^j_{p_{i}^{l}}-W_1\big)\cdot H_1dx\\
&-\int_{\R^2}\nabla_{i}{W}^j_{p_{i}^{l}}\cdot\big(\Delta_j\delta_{\mn\times\mh\times\mn}-H_1\big)dx\\
\triangleq& -\int_{\R^2}W_1\cdot H_1dx+B_1+B_6.
\end{align*}
Due to Lemma \ref{lem:W1 H1 positivity},
\begin{align*}
\f14\int_{\R^2}W_1\cdot H_1dx
\ge\frac{3}{4}a^2\|\Delta\Delta_j\delta_{\mn}\|_{2}^2-B_1,
\end{align*}
which along with (\ref{eq:bar w}) gives the lemma.\endproof\vspace{0.1cm}

Now we follow the same route as in \cite{WangWZ2013} and begin with the estimates of $B_i$. \vspace{0.1cm}

$\bullet$ \underline{Estimate of $B_1$.}

By (\ref{eq:nabla W j}) and the definition of $W_1$ and $H_1$, we have
\begin{align*}
B_5\le& C\|\na \mn^2\|_4\|\Delta_j\na\delta_\mn\|_4\|\Delta_j\na^2\delta_\mn\|_2
\le C\|\nabla\mn^2\|_4\|\Delta_j\nabla\delta_{\mn}\|_2^{\f12}\|\Delta_j\nabla^2\delta_{\mn}\|_2^{\f32}\\
\le& C\|\nabla\mn^2\|_{4}^{4}\|\Delta_j\nabla\delta_{\mn}\|_2^2+\epsilon2^{4j}\|\Delta_j\delta_{\mn}\|_2^2\\
\le& \epsilon 2^{4j}\|\Delta_j\delta_{\mn}\|_{2}^2+C2^{2js}\bar{h}(t)W(t).
\end{align*}

$\bullet$ \underline{Estimate of $B_2$.}
Recall that
\ben\label{eq:bar wn}
{W}_{\mn_l^2}^{j}=2(k_3-k_2)(\mn^2\cdot\mcurl\Delta_{j}\delta_{\mn})\mcurl\Delta_{j}\delta_{\mn},
\een
then Lemma \ref{lem:Berstein} yields that
\begin{align*}
B_1\leq& C2^{2j}\|\Delta_j\delta_{\mn}\|_{2}\|\partial_t{\mn}^2\|_{2}\|\Delta_j\delta_{\mn}\|_{\infty}\le
C2^{3j}\|\Delta_j\delta_{\mn}\|_{2}^2\|\partial_t{\mn}^2\|_{2}\nonumber\\
\leq& \epsilon 2^{4j}\|\Delta_j\delta_{\mn}\|_{2}^2+C2^{2j}\|\partial_t{\mn}^2\|_{2}^2\|\Delta_j\delta_{\mn}\|_{2}^2\\
\le&  \epsilon 2^{4j}\|\Delta_j\delta_{\mn}\|_{2}^2+C2^{2js}\bar{h}(t)W(t).
\end{align*}

$\bullet$ \underline{Estimate of $B_6$, $B_3,B_4,B_5$.}

By (\ref{eq:nabla W j}) and Lemma \ref{lem:Berstein}, for $j\geq 0$ we have
\begin{align}
\|\nabla_{i}{W}^j_{p_{i}^{l}}\|_{2}\leq  C\big(\|\nabla^2\Delta_j\delta_{\mn}\|_{2}+\||\nabla\mn^2|\Delta_j\delta_{\nabla\mn}\|_2\big)
\leq C2^{2j}\|\Delta_j\delta_{\mn}\|_{2}.\label{eq:dW-est}
\end{align}
Denote $\cB$  the following form
\beno
\sum_{i=1,2}\Delta_j(\cM\nabla^2\mn^i\delta_{\mn})+\Delta_j(\cM\nabla\mn^i\nabla\delta_{\mn})+\sum_{i,k=1,2}\Delta_j(\cM\nabla\mn^i\nabla\mn^k\delta_{\mn})
+[\Delta_j,\cM]\nabla^2\delta_{\mn}.
\eeno
Then by (\ref{eq:hn-d}) and (\ref{eq:dW-est}), we have
\beno
B_6\le C2^{2j}\|\Delta_j\delta_\mn\|_2\|\cB\|_2,
\eeno
where $\|\cB\|_2$ is bounded by
\beno
\sum_{i=1,2}\|\Delta_j(\cM\nabla^2\mn^i\delta_{\mn})\|_2+\|\Delta_j(\cM\nabla\mn^i\nabla\delta_{\mn})\|_2+\sum_{i,k=1,2}\|\Delta_j(\cM\nabla\mn^i\nabla\mn^k\delta_{\mn})\|_2
+\|[\Delta_j,\cM]\nabla^2\delta_{\mn}\|_2.
\eeno
Then it follows from Lemma \ref{lem:product3-1}--Lemma \ref{commutator} that
\beno
B_6\le \epsilon \sum_{l=j-9}^{j+9}2^{4l}\|\Delta_{l}\delta_{\mn}\|_{2}^2
+C2^{2js}\bar{h}(t)W(t).
\eeno
Moreover,
\beno
|B_3|+|B_4|+|B_5|\le C2^{2j}\|\Delta_j\delta_\mn\|_2\|\cB\|_2,
\eeno
and
\beno
|B_3|+|B_4|+|B_5|+|B_6|\le 4\epsilon \sum_{l=j-9}^{j+9}2^{4l}\|\Delta_{l}\delta_{\mn}\|_{2}^2
+C2^{2js}\bar{h}(t)W(t).
\eeno

Thus, Proposition \ref{prop:W} follows from Lemma \ref{lem:W-evolution} and the estimates for $B_i$.\endproof


\appendix

\setcounter{equation}{0}
\section{Local well-posedness results in $\R^d$ with $d=2,3$}
The symbol $\langle\cdot,\cdot\rangle$ denotes the integral in $\R^d$ with $d=2,3$. Moreover, ${\mathcal P}(\cdot,\cdots,\cdot)$ denotes \textcolor[rgb]{0.00,0.00,0.00}{a} polynomial depending on the \textcolor[rgb]{0.00,0.00,0.00}{arguments}  in the \textcolor[rgb]{0.00,0.00,0.00}{parentheses} whose order, for example, is less than $10$.

In this section, we are aimed to prove the local existence and  blow-up criterion  for strong solutions of the system \eqref{eq: Landau-Lifshitz with Oseen-Frank} in $\R^d$ with $d=2,3$. Firstly, we use the classical  Friedrich's method to construct the approximate solutions of (\ref{eq: Landau-Lifshitz with Oseen-Frank}) as in \cite{WZZ2013,WangW2014}. The main difficulty lies in the Schr\"{o}dinger term $\mn\times\mh$, which can't be controlled by the term $\mn\times(\mn\times\mh)$ when $|\mn|\neq 1$. Hence, we introduce an equivalent system of   (\ref{eq: Landau-Lifshitz with Oseen-Frank}) as follows
\ben\label{eq:L-L Oseen-Frank=}
\partial_t\mn= \alpha\mn\times(\mh\times\mn) +\beta\mn\times[(\mn\times\mh)\times\mn],
\een
Secondly, blow-up criterion  is similar to \cite{WangW2014}. We'll use a better representation formula of $\mh\cdot\mn$ and the vertical property of $\mn\times\mh$ with respect to $\mn$.

Our main theorem states as follows.

\begin{Theorem}\label{thm:local existence}
Let $s\ge 2$  be an integer, and
the initial data $\nabla\mn_0\in H^{2s}(\R^d)$ for $d=2,d=3$.
Then there exist $T>0$ and a solution $\mn$ of the system (\ref{eq: Landau-Lifshitz with Oseen-Frank}) such that
\beno
\nabla\mn\in C\big([0,T^*);H^{2s}(\R^d)\big).
\eeno
Moreover, if $T^{*}$ is the maximal existence time of the solution, then $T^{*}<+\infty$ implies that
$$
\int_{0}^{T^{*}}\|\nabla\mn(t)\|_{L^{\infty}}^{2}dt =+\infty.
$$
\end{Theorem}


The following lemma will be frequently used for the commutator; for example see \cite{bcd}.
\begin{Lemma}\label{lem:galiado}
For $\alpha,\beta\in N^3$ or $N^2$, it holds that
$$
\|D^{\alpha}(fg)\|_{L^2}\le C\sum_{|\gamma|=|\alpha|}\big(\|f\|_{L^{\infty}}\|D^{\gamma}g\|_{L^2}+\|g\|_{L^{\infty}}\|D^{\gamma}f\|_{L^2}\big),
$$
$$
\|[D^{\alpha},f]D^{\beta}g\|_{L^2}\le
C\left(\sum_{|\gamma|=|\alpha|+|\beta|}\|D^{\gamma}f\|_{L^2}\|g\|_{L^{\infty}}+\sum_{|\gamma|=|\alpha|+|\beta|-1}\|\nabla
f\|_{L^{\infty}}\|D^{\gamma}g\|_{L^2}\right).
$$
\end{Lemma}

Let $a,k_1,k_2,k_3$ be the parameters of $\mh$, then we have the following inequality.
\begin{Lemma}\label{lem:positive of h}For any vector $f\in L^2(\R^d)$, there holds
\begin{align*}
 (k_2-a)\|f\|_2^2+(k_3-k_2)\|\mn\cdot f\|_2^2\geq 0.
\end{align*}
\end{Lemma}

In fact, on one hand
\begin{itemize}
  \item if $a=k_1$, then either $k_3\geq k_2$ or $|k_3-k_2|\leq |k_2-a|$;
  \item if $a=k_2$, then $k_3\geq k_2$;
  \item if $a=k_3$, then $|k_3-k_2|=|k_2-a|$,
\end{itemize}
on the other hand,
$|{\bf n}\times\mcurl{\bf n}|^2+|{\bf n}\cdot(\mcurl{\bf n})|^2=|\mcurl{\bf n}|^2$ implies the above inequality.

{\bf Proof of Theorem \ref{thm:local existence}:} It's divided into three steps.

{\bf Step 1. Construction of the approximated solutions:}
Let $b\in S^2$ be a constant vector, $\mn_0:\R^d\to S^2$ such that $\mn_0-b\in H^{k}(\R^d)$ with $k>0$.
Let
$$
{\cal J}_{\epsilon}f={\cal F}^{-1}(\phi( {\frac{\xi}{\epsilon}}){\cal F}f),
$$
where ${\cal F}(f)(\xi)=\int_{\R^d}f(x)e^{-ix\xi}dx$ is usual Fourier transform and $\phi(\xi)$ is a smooth cut-off function with $\phi=1$ in $B_1$ and
$\phi=0$ outside of $B_2$. We construct the approximate system of (\ref{eq:L-L Oseen-Frank=}),
\begin{eqnarray}\label{eq:approximate system of LL}
\left\{\begin{array}{lll}\partial_t\mn_{\epsilon}=\alpha{\cal J}_{\epsilon}\left({\cal J}_{\epsilon}\mn_{\epsilon}\times({\cal
J}_{\epsilon}h_{\epsilon}\times{\cal J}_{\epsilon}\mn_{\epsilon})\right)+\beta{\cal J}_{\epsilon}\left({\cal J}_{\epsilon}\mn_{\epsilon}\times\left(({\cal
J}_{\epsilon}\mn_{\epsilon}\times{\cal J}_{\epsilon}h_{\epsilon})\times{\cal J}_{\epsilon}\mn_{\epsilon}\right)\right),\\
\mn_{\epsilon}|_{t=0}={\cal J}_{\epsilon}\mn_0,\end{array}\right.
\end{eqnarray}
where
\begin{eqnarray*}{\cal J}_{\epsilon}h_{\epsilon}&=&2a\Delta{\cal J}_{\epsilon}\mne+2(k_1-a)\nabla{\rm div}{\cal
J}_{\epsilon}\mne-2(k_2-a)\mcurl(\je\mne \times(\mcurl\je\mne\times\je\mne))\\
&&-2(k_3-a)\mcurl((\je\mne\cdot\mcurl\je\mne)\je\mne) -2(k_3-k_2)(\mcurl\je\mne\cdot\je\mne) \cdot\mcurl\je\mne.
\end{eqnarray*}

By the Cauchy-Lipschitz theorem (for example, see \cite{Ar}), we know that there exists a strictly maximal time $T_{\epsilon}$ and
a unique solution $\mn_{\epsilon}-\mn_0\in C([0,T_{\epsilon});H^{k}(\R^d))$ for any $k> 0$.

{ \bf Step 2. Uniform energy estimates:}
We consider the evolution of the following energy norm
\begin{eqnarray*}
E_{s}(\mne)&= &\|\mne-\mn_0\|_{2}^{2} +\int_{\R^d}W(\mne,\nabla\mne)dx +a\|\Delta^s\nabla\mne\|_2^2+(k_1-a)\|\Delta^s
div\mne\|_2^2\nonumber\\&& +(k_2-a)\|\Delta^s\mcurl\mne\times\je\mne\|_2^2+(k_3-a)\|\Delta^s\mcurl\mne\cdot\je\mne\|_2^2,
\end{eqnarray*}
and it's sufficient to prove that
\ben\label{eq:uniform energy estimates}
\frac{d}{dt}E_{s}(\mne)\le  C{\mathcal P}(\|\mn_{\epsilon}\|_{L^{\infty}},\|\nabla\mn_{\epsilon}\|_{L^{\infty}},
\|\nabla^2\mn_{\epsilon}\|_{L^{\infty}})E_{s}(\mne)\le {\cal F}(E_{s}(\mne)),
\een
where we used the embedding equality with $s\geq 2$ and ${\cal F}$ is an increasing function with ${\cal F}(0)=0$. Indeed, it means that there exists a $T>0$ depending only on $E_{s}(\mn_0)$ such that for all
$t\in[0,\min(T,T_{\epsilon})]$,
\begin{eqnarray*}
E_{s}(\mne)\le 2E_{s}(\mn_{0}),
\end{eqnarray*}
which implies that $T_{\epsilon}\ge T$ by a continuous argument. Then the uniform estimates for the solutions $\mne$ on $[0,T]$ hold which yield that
there exists a local solution $\mn$ of (\ref{eq:L-L Oseen-Frank=}) by the standard compactness arguments. Also, if $|\mn_0|=1$, multiply
$\cdot\mn$ on both sides of (\ref{eq:L-L Oseen-Frank=}) and we can obtain $|\mn|=1$.

Next, we come to prove the estimate (\ref{eq:uniform energy estimates}).

{\bf  2.1. Lower order terms:} In fact, using the equation (\ref{eq:approximate system of LL}) we have
\begin{align*}
\frac{1}{2}\frac{d}{dt}\|\mn_\epsilon-\mn_0\|_2^2 =&\langle\partial_t\mn_{\epsilon},\mn_{\epsilon}-\mn_0\rangle\nonumber\\
\le &C(1+\|\mne\|_{L^\infty}+\|\nabla\mne\|_{L^\infty})^4(\|\nabla\mne\|_{2}+\|\Delta\mne\|_2)\|\mne-\mn_0\|_2\nonumber\\
\le& C(1+\|\mn_{\epsilon}\|_{L^\infty}+\|\nabla\mne\|_{L^\infty})^4E_{s}(\mn_{\epsilon})
\end{align*}
and on the other hand
\begin{align*}
&\frac{d}{dt}\int_{\R^d}W(\mne,\nabla\mne)(\cdot,t)dx\\
=&\int_{\R^d}\left(W_{n^l}-\nabla_iW_{p_{i}^{l}}\right) (\mne,\nabla\mne)\nonumber\\
& \cdot\je\left(\alpha(\je\mne\times(\je h_{\epsilon}\times\je\mne))+\beta\left({\cal J}_{\epsilon}\mn_{\epsilon}\times\left(({\cal
J}_{\epsilon}\mn_{\epsilon}\times{\cal J}_{\epsilon}h_{\epsilon})\times{\cal J}_{\epsilon}\mn_{\epsilon}\right)\right)\right)dx\nonumber\\
\le & C{\mathcal P}\big(\|\mn_{\epsilon}\|_{L^{\infty}},\|\nabla\mn_{\epsilon}\|_{L^{\infty}}\big)
\|\nabla\mn_{\epsilon}\|_{H^{1}}^2,
\end{align*}
which are the required estimates.

{\bf 2.2. Higher order term:} Direct calculation shows that
\begin{align*}
\frac{1}{2}\frac{d}{dt}\langle\nabla\Delta^s\mne,\nabla\Delta^s\mne\rangle=&-\alpha\langle\Delta^s(\je\mne\times(\je
h_{\epsilon}\times\je\mne)),\je\Delta^{s+1}\mne\rangle\nonumber\\
&-\beta\langle\Delta^s({\cal J}_{\epsilon}\mn_{\epsilon}\times\left(({\cal
J}_{\epsilon}\mn_{\epsilon}\times{\cal J}_{\epsilon}h_{\epsilon})\times{\cal J}_{\epsilon}\mn_{\epsilon}\right)),\je\Delta^{s+1}\mne\rangle\nonumber\\
:=&I_1+I_2,
\end{align*}
\begin{align*}
\frac{1}{2}\frac{d}{dt}\langle div\Delta^s\mne, div\Delta^s\mne\rangle=&-\alpha\langle\Delta^s(\je\mne\times(\je h_{\epsilon}\times\je\mne)),\je\nabla
div\Delta^{s}\mne\rangle\nonumber\\
&-\beta\langle\Delta^s({\cal J}_{\epsilon}\mn_{\epsilon}\times\left(({\cal
J}_{\epsilon}\mn_{\epsilon}\times{\cal J}_{\epsilon}h_{\epsilon})\times{\cal J}_{\epsilon}\mn_{\epsilon}\right)),\je\nabla div\Delta^{s}\mne\rangle\nonumber\\
:=&I_1'+I_2',
\end{align*}
\begin{align*}
&\frac{1}{2}\frac{d}{dt}\langle\je\mne\times\Delta^s\mcurl\mne,\je\mne\times\Delta^s\mcurl\mne\rangle\nonumber\\
=&\alpha\langle\je\mne\times\Delta^s\mcurl\je\left(\je\mne\times(\je h_{\epsilon}\times\je\mne)\right),\je\mne\times\Delta^s\mcurl\mne\rangle\nonumber\\
&+\beta\langle\je\mne\times\Delta^s\mcurl\je\left({\cal J}_{\epsilon}\mn_{\epsilon}\times\left(({\cal
J}_{\epsilon}\mn_{\epsilon}\times{\cal J}_{\epsilon}h_{\epsilon})\times{\cal J}_{\epsilon}\mn_{\epsilon}\right)\right),\je\mne\times\Delta^s\mcurl\mne\rangle\nonumber\\
&+\alpha\langle\je^2\left(\je\mne\times(\je h_{\epsilon}\times\je\mne)\right)\times\Delta^s\mcurl\mne,\je\mne\times\Delta^s\mcurl\mne\rangle\nonumber\\
&+\beta\langle\je^2\left({\cal J}_{\epsilon}\mn_{\epsilon}\times\left(({\cal
J}_{\epsilon}\mn_{\epsilon}\times{\cal J}_{\epsilon}h_{\epsilon})\times{\cal J}_{\epsilon}\mn_{\epsilon}\right)\right)\times\Delta^s\mcurl\mne,\je\mne\times\Delta^s\mcurl\mne\rangle\nonumber\\
:=& I_1''+I_2''+I_3''+I_4'',
\end{align*}
and
\begin{align*}
&\frac{1}{2}\frac{d}{dt}\langle\je\mne\cdot\Delta^s\mcurl\mne,\je\mne\cdot\Delta^s\mcurl\mne\rangle\nonumber\\
=&\alpha\langle\je\mne\cdot\Delta^s\mcurl\je\left(\je\mne\times(\je h_{\epsilon}\times\je\mne)\right),\je\mne\cdot\Delta^s\mcurl\mne\rangle\nonumber\\
&+\beta\langle\je\mne\cdot\Delta^s\mcurl\je\left({\cal J}_{\epsilon}\mn_{\epsilon}\times\left(({\cal
J}_{\epsilon}\mn_{\epsilon}\times{\cal J}_{\epsilon}h_{\epsilon})\times{\cal J}_{\epsilon}\mn_{\epsilon}\right)\right),\je\mne\cdot\Delta^s\mcurl\mne\rangle\nonumber\\
&+\alpha\langle\je^2\left(\je\mne\times(\je h_{\epsilon}\times\je\mne)\right)\cdot\Delta^s\mcurl\mne,\je\mne\cdot\Delta^s\mcurl\mne\rangle\nonumber\\
&+\beta\langle\je^2\left({\cal J}_{\epsilon}\mn_{\epsilon}\times\left(({\cal
J}_{\epsilon}\mn_{\epsilon}\times{\cal J}_{\epsilon}h_{\epsilon})\times{\cal J}_{\epsilon}\mn_{\epsilon}\right)\right)\cdot\Delta^s\mcurl\mne,\je\mne\cdot\Delta^s\mcurl\mne\rangle\nonumber\\
:=&
I_1'''+I_2'''+I_3'''+I_4'''.
\end{align*}
Then we have
\begin{align*}
&2(k_2-a)I_{3}''+2(k_2-a)I_4''+2(k_3-a)I_3'''+2(k_3-a)I_4'''\nonumber\\
\le& C{\mathcal P}(\|\mn_{\epsilon}\|_{L^{\infty}},\|\nabla\mn_{\epsilon}\|_{L^{\infty}},
\|\nabla^2\mn_{\epsilon}\|_{L^{\infty}})\|\nabla\mn_{\epsilon}\|_{H^{2s}}^2.\nonumber
\end{align*}

By the formula of $\je\mne$ and commutator estimates in Lemma \ref{lem:galiado}, we get
\begin{align}
&\|\Delta^s(\je\mne\times(\je h_{\epsilon}\times\je\mne))\|_{L^2}\nonumber\\
\le &\|[\Delta^s,\je\mne\times](\je h_{\epsilon}\times\je\mne)\|_{L^2}+\|\je\mne\times\Delta^s(\je h_{\epsilon}\times\je\mne)\|_{L^2}\nonumber\\
\le& {\cal
P}(\|\mne\|_{L^{\infty}},\|\nabla\mne\|_{L^{\infty}},\|\nabla^2\mne\|_{L^{\infty}})\left(\|\nabla\mne\|_{H^{2s}}+\|\Delta^s\je h_{\epsilon}\times\je\mne\|_2\right),
\end{align}
Recall the commutator estimates of $[\je,f]$ in \cite{WangW2014},
\begin{align*}
\|[{\cal J}_{\epsilon},f]\nabla g\|_{L^p}
\leq& C(1+\|\nabla f\|_{L^{\infty}})\|g\|_{L^p},
\end{align*}
therefore we have
\begin{align*}
&2aI_{1}+2(k_1-a)I_1'+2(k_2-a)I_1''+2(k_3-a)I_1'''\nonumber\\
=&\alpha\langle\Delta^s(\je\mne\times(\je h_{\epsilon}\times\je\mne)),-2a\je\Delta^{s+1}\mne-2(k_1-a)\je\nabla div\Delta^s\mne\nonumber\\
&~~~~~~~~~~~~~~~~~~~~~~~~~~~~~~~~~~~~~~+2(k_2-a)\je\mcurl\left((\je\mne\times\Delta^s\mcurl\mne)\times\je\mne\right)\nonumber\\
&~~~~~~~~~~~~~~~~~~~~~~~~~~~~~~~~~~~~~~+2(k_3-a)\je\mcurl\left((\je\mne\cdot\Delta^s\mcurl\mne)\cdot\je\mne\right)\rangle\nonumber\\
\le &-\alpha\langle\Delta^s(\je\mne\times(\je h_{\epsilon}\times\je\mne)),\Delta^s\je h_{\epsilon}\rangle\nonumber\\&
+C(\delta){\mathcal P}(\|\mn_{\epsilon}\|_{L^{\infty}},\|\nabla\mn_{\epsilon}\|_{L^{\infty}},
\|\nabla^2\mn_{\epsilon}\|_{L^{\infty}}) \|\nabla\mn_{\epsilon}\|_{H^{2s}}^2+\delta\|\Delta^s\je h_{\epsilon}\times\je\mne\|_{L^2}^{2}\nonumber\\
\le & -\frac{\alpha}{2}\langle\Delta^s \je h_{\epsilon}\times\je\mne,\Delta^s \je h_{\epsilon}\times\je\mne\rangle+C{\mathcal
P}(\|\mn_{\epsilon}\|_{L^{\infty}},\|\nabla\mn_{\epsilon}\|_{L^{\infty}},
\|\nabla^2\mn_{\epsilon}\|_{L^{\infty}}) \|\nabla\mn_{\epsilon}\|_{H^{2s}}^2
\end{align*}
Similarly,
\begin{align*}
&2aI_{2}+2(k_1-a)I_2'+2(k_2-a)I_2''+2(k_3-a)I_2'''\nonumber\\
\le &-\beta\langle\Delta^s \left({\cal J}_{\epsilon}\mn_{\epsilon}\times\left(({\cal
J}_{\epsilon}\mn_{\epsilon}\times{\cal J}_{\epsilon}h_{\epsilon})\times{\cal J}_{\epsilon}\mn_{\epsilon}\right)\right),\Delta^s \je h_{\epsilon}\rangle\nonumber\\
&+C(\delta){\mathcal
P}(\|\mn_{\epsilon}\|_{L^{\infty}},\|\nabla\mn_{\epsilon}\|_{L^{\infty}},
\|\nabla^2\mn_{\epsilon}\|_{L^{\infty}}) \|\nabla\mn_{\epsilon}\|_{H^{2s}}^2+\delta\|\Delta^s\je h_{\epsilon}\times\je\mne\|_{L^2}^{2}\nonumber\\
\le& C(\delta){\mathcal
P}(\|\mn_{\epsilon}\|_{L^{\infty}},\|\nabla\mn_{\epsilon}\|_{L^{\infty}},
\|\nabla^2\mn_{\epsilon}\|_{L^{\infty}}) \|\nabla\mn_{\epsilon}\|_{H^{2s}}^2+2\delta\|\Delta^s\je h_{\epsilon}\times\je\mne\|_{L^2}^{2},
\end{align*}
where we used the relation
\beno
\langle \left({\cal J}_{\epsilon}\mn_{\epsilon}\times\left(({\cal
J}_{\epsilon}\mn_{\epsilon}\times\Delta^s{\cal J}_{\epsilon}h_{\epsilon})\times{\cal J}_{\epsilon}\mn_{\epsilon}\right)\right),\Delta^s \je h_{\epsilon}\rangle=0,
\eeno
and thus
\begin{align*}
&\langle\Delta^s \left({\cal J}_{\epsilon}\mn_{\epsilon}\times\left(({\cal
J}_{\epsilon}\mn_{\epsilon}\times{\cal J}_{\epsilon}h_{\epsilon})\times{\cal J}_{\epsilon}\mn_{\epsilon}\right)\right),\Delta^s \je h_{\epsilon}\rangle\\
=&\langle[\Delta^s,\je\mne\times(\je\mne\times)](\je h_{\epsilon}\times\je\mne),\Delta^s\je h_{\epsilon}\rangle\nonumber\\
&-\langle[\Delta^s,(\je \mne\times)]\je h_{\epsilon},(\Delta^s\je h_{\epsilon}\times\je\mne)\times\je\mne\rangle\nonumber\\
\le&\langle[\Delta^s,\je\mne\times(\je\mne\times)](\nabla\je h_{\epsilon}\times\je\mne),\nabla\Delta^{s-1}\je h_{\epsilon}\rangle\nonumber\\
&+C(\delta){\mathcal
P}(\|\mn_{\epsilon}\|_{L^{\infty}},\|\nabla\mn_{\epsilon}\|_{L^{\infty}},
\|\nabla^2\mn_{\epsilon}\|_{L^{\infty}}) \|\nabla\mn_{\epsilon}\|_{H^{2s}}^2+\delta\|\Delta^s\je h_{\epsilon}\times\je\mne\|_{L^2}^{2}\\
\leq&C(\delta){\mathcal
P}(\|\mn_{\epsilon}\|_{L^{\infty}},\|\nabla\mn_{\epsilon}\|_{L^{\infty}},
\|\nabla^2\mn_{\epsilon}\|_{L^{\infty}}) \|\nabla\mn_{\epsilon}\|_{H^{2s}}^2+2\delta\|\Delta^s\je h_{\epsilon}\times\je\mne\|_{L^2}^{2}.
\end{align*}
Combining the above estimates, the inequality (\ref{eq:uniform energy estimates})  is satisfied by choosing $\delta$ is sufficiently small, thus the proof of the local existence is complete.

{ \bf Step 3. Blow-up criterion.}

Let  $T^{*}<\infty$ be the maximal existence time of the solution. Then it is sufficient to prove that
\ben\label{eq:blow-up estimate}
\frac{d}{dt}E_{s}(\mn)\le C(1+\|\nabla\mn\|_{L^{\infty}}^{2})E_{s}(\mn),
\een
where
\begin{align*}
E_{s}(\mn)=&\|\mn-\mn_0\|_{L^2}^2+\int_{\R^d}W(\mn,\nabla\mn)dx+a\|\Delta^s\nabla\mn\|_{L^2}^2
+(k_1-a)\|\Delta^s{\rm div}\mn\|_{L^2}^2\\
&+(k_2-a)\|\mn\times\Delta^s(\nabla\times\mn)\|_{L^2}^{2}+(k_3-a)\|\mn\cdot\Delta^s(\nabla\times\mn)\|_{L^2}^{2}.
\end{align*}
The proof of (\ref{eq:blow-up estimate}) is more subtle with respect to the existence, since we can't use the bound of $\|\nabla^2\mn\|_{\infty}$. However, at this time we have $|\mn|=1$, and
$\mn\cdot\Delta\mn=-|\nabla\mn|^2$.

{\bf  3.1. Lower order terms:} It is easy to see that
\begin{align*}
\frac{d}{dt}\|\mn-\mn_0\|_2^2=&\langle\partial_t\mn,\mn-\mn_0\rangle\nonumber\\
=&2\langle\alpha\mn\times(\mh\times\mn)+\beta\mn\times\mh,\mn-\mn_0\rangle\nonumber\\
\le &C(\|\Delta\mn\|_2+\||\nabla\mn|^2\|_2)\|\mn-\mn_0\|_2\le CE_{s}(\mn),
\end{align*}
and
\beno
\frac{d}{dt}\int_{\R^d}W(\mn,\nabla\mn)(\cdot,t)dt
=\int_{\R^d}\left(W_{n^l}-\nabla_iW_{p_{i}^{l}}\right)\partial_{t}n^ldx
=-\alpha\int_{\R^d}|\mn\times\mh|^2dx.
\eeno

{\bf  Step 3.2. Higher order term:} Direct calculation shows that
\begin{align*}
&\frac{1}{2}\frac{d}{dt}\langle\nabla\Delta^s\mn,\nabla\Delta^s\mn\rangle\\
=&-\alpha\langle\Delta^s(\mn\times(
\mh\times\mn)),\Delta^{s+1}\mn\rangle-\beta\langle\Delta^s(\mn\times\mh),\Delta^{s+1}\mn\rangle:=I_1+I_2,
\end{align*}
\begin{align*}
&\frac{1}{2}\frac{d}{dt}\langle div\Delta^s\mn, div\Delta^s\mn\rangle\\
=&-\alpha\langle\Delta^s(\mn\times(\mh\times\mn)),\nabla
div\Delta^{s}\mn\rangle-\beta\langle\Delta^s(\mn\times\mh),\nabla div\Delta^{s}\mn\rangle:=I_1'+I_2',
\end{align*}
\begin{align*}
&\frac{1}{2}\frac{d}{dt}\langle\mn\times\Delta^s\mcurl\mn,\mn\times\Delta^s\mcurl\mn\rangle\\
=&\alpha\langle\mn\times\Delta^s\mcurl\left(\mn\times(\mh\times\mn)\right),\mn\times\Delta^s\mcurl\mn\rangle\\
&+\beta\langle\mn\times\Delta^s\mcurl\left(\mn\times\mh\right),\mn\times\Delta^s\mcurl\mn\rangle\\
&+\alpha\langle\left(\mn\times(\mh\times\mn)\right)\times\Delta^s\mcurl\mn,\mn\times\Delta^s\mcurl\mn\rangle\\
&+\beta\langle\left(\mn\times\mh\right)\times\Delta^s\mcurl\mn,\mn\times\Delta^s\mcurl\mn\rangle\\
:=& I_1''+I_2''+I_3''+I_4'',
\end{align*}
and
\begin{align*}
&\frac{1}{2}\frac{d}{dt}\langle\mn\cdot\Delta^s\mcurl\mn,\mn\cdot\Delta^s\mcurl\mn\rangle\\
=&\alpha\langle\mn\cdot\Delta^s\mcurl\left(\mn\times(\mh\times\mn)\right),\mn\cdot\Delta^s\mcurl\mn\rangle\\
&+\beta\langle\mn\cdot\Delta^s\mcurl\left(\mn\times\mh\right),\mn\cdot\Delta^s\mcurl\mn\rangle\\
&+\alpha\langle\left(\mn\times(\mh\times\mn)\right)\times\Delta^s\mcurl\mn,\mn\cdot\Delta^s\mcurl\mn\rangle\\
&+\beta\langle\left(\mn\times\mh\right)\cdot\Delta^s\mcurl\mn,\mn\cdot\Delta^s\mcurl\mn\rangle\\
:=&
I_1'''+I_2'''+I_3'''+I_4'''.
\end{align*}
For the terms $I_1,I_1',I_1'',I_1'''$, we have
\begin{align}
&2aI_1+2(k_1-a)I_1'+2(k_2-a)I_1''+2(k_3-a)I_1'''\nonumber\\
=&\alpha\langle\Delta^s(\mn\times(\mh\times\mn)),-2a\Delta^{s+1}\mn-2(k_1-a)\nabla div\Delta^s\mn+2(k_2-a)\mcurl\mcurl\Delta^s\mn\rangle\nonumber\\
&+\alpha\langle\Delta^s(\mn\times(\mh\times\mn)),2(k_3-k_2)[(\mn\cdot\Delta^s\mcurl\mn)\mcurl\mn+\nabla(\mn\cdot\Delta^s\mcurl\mn)\times\mn]\rangle\nonumber\\
=&-\alpha\langle\Delta^s(\mn\times(\mh\times\mn)),\Delta^s\nabla_{\alpha}W_{p_{\alpha}^{l}}\rangle\nonumber\\
&+2(k_3-k_2)\alpha\langle\Delta^s(\mn\times(\mh\times\mn)),(\mn\cdot\Delta^s\mcurl\mn)\mcurl\mn-\Delta^s\left((\mn\cdot\mcurl\mn)\mcurl\mn\right)\rangle\nonumber\\
&+2(k_3-k_2)\alpha\langle\Delta^s(\mn\times(\mh\times\mn)),\nabla(\mn\cdot\Delta^s\mcurl\mn)\times\mn-\nabla\Delta^s(\mn\cdot\mcurl\mn)\times\mn\rangle\nonumber\\
&+2(k_3-k_2)\alpha\langle\nabla\Delta^s(\mn\cdot\mcurl\mn)\times\mn-\Delta^s\left(\nabla(\mn\cdot\mcurl\mn)\times\mn\right)\rangle,
\end{align}
where we have used the following relation, for a function $f$ and a vector field $u$, there holds
$$
\mcurl(fu)=f\mcurl u+\nabla f\times u.
$$

We will use the following Gagliardo-Sobolev inequality on $\R^d$ (for example, see \cite{Ad}). Let $\tau\in N$, and $\tau\ge 2s-1$, then for $1\le j\le [\tau/2]$, $[\tau/2]+1\le k\le \tau$, and $f\in H^{\tau+1}(\R^d)$, we have
\begin{eqnarray*}
\|\nabla^j f\|_{L^{\infty}}\le C\|\nabla f\|_{H^{\tau}}^{\frac{j}{\tau+1-d/2}}\|f\|_{L^{\infty}}^{1-\frac{j}{\tau+1-d/2}},\\
\|\nabla^k f\|_{L^2}\le C\|\nabla f\|_{H^{\tau}}^{\frac{k-d/2}{\tau+1-d/2}}\|f\|_{L^{\infty}}^{1-\frac{k-d/2}{\tau+1-d/2}}.
\end{eqnarray*}
Hence, for $\tau\ge 2s-1$ with $s\ge 2$, the following inequalities hold,
\begin{eqnarray}\label{gs}
\|\nabla^{\tau+1}\mn\|_{L^2}\|\nabla\mn\|_{L^{\infty}}+\|\nabla^{\tau}\mn\|_{L^2}\|\nabla^2\mn\|_{L^{\infty}}
+\|\nabla^{\tau}\mn\|_{L^2}\|\nabla\mn\|_{L^{\infty}}^2\le C \|\nabla\mn\|_{H^{\tau+1}},\nonumber\\
\left(\|\nabla^2\mn\|_{L^{\infty}}+\|\nabla\mn\|_{L^{\infty}}^{2}\right)\|\nabla^{\tau}\mn\|_{L^2}\le C\|\nabla\mn\|_{L^{\infty}}\|\nabla^{\tau+1}\mn\|_{L^2}.
\end{eqnarray}
By Lemma \ref{lem:galiado} and Gagliardo-Sobolev inequality (\ref{gs}), we have
\begin{align}\label{ii5'''}
&2aI_1+2(k_1-a)I_1'+2(k_2-a)I_1''+2(k_3-a)I_1'''\nonumber\\
\le& -\alpha\langle\Delta^s(\mn\times(\mh\times\mn)),\Delta^s\nabla_{\alpha}W_{p_{\alpha}^{l}}\rangle\nonumber\\
&+C\big(\|\Delta^s\mn\|_{L^2}
+\|\Delta^s\mn\|_{L^2}(\|\nabla\mn\|_{L^{\infty}}^2+\|\nabla^2\mn\|_{L^{\infty}})+\|\Delta^s\mh\times\mn\|_{L^2}\big)\nonumber\\
&\cdot
(\|\Delta^s\mn\|_{L^2}\|\nabla\mn\|_{L^{\infty}}^2+\|\nabla\mn\|_{L^{\infty}}\|\Delta^s\nabla\mn\|_{L^2}+\|\Delta^s\mn\|_{L^2}\|\nabla^2\mn\|_{L^{\infty}})\nonumber\\
\le&
-\alpha\langle\Delta^s(\mn\times(\mh\times\mn)),\Delta^s\nabla_{\alpha}W_{p_{\alpha}^{l}}\rangle+C_{\delta}\|\nabla\mn\|_{L^{\infty}}^2\|\Delta^s\nabla\mn\|_{L^2}^2+\delta\|\Delta^{s+1}\mn\|_{L^2}^2.
\end{align}

Note that
\begin{align}\label{i}
&-\alpha\langle\Delta^s(\mn\times(\mh\times\mn)),\Delta^s\nabla_{\alpha}W_{p_{\alpha}^{l}}\rangle\nonumber\\
=&-\alpha\big\langle\Delta^s\left(\mn\times\left((\mh-\nabla_{\alpha}W_{p_{\alpha}^{k}})\times\mn\right)\right),\Delta^s\nabla_{\alpha}W_{p_{\alpha}^{l}}\big\rangle\nonumber\\
&-\alpha\langle\Delta^{s}\left(\mn\times(\nabla_{\alpha}W_{p_{\alpha}^{k}}\times\mn)\right),
\Delta^s(2a\Delta \mn)\rangle\nonumber\\
&-\alpha\langle\Delta^{s}
\left(\mn\times((\nabla_{\alpha}W_{p_{\alpha}^{k}}-2a\Delta \mn)\times\mn)\right),
\Delta^s(\nabla_{\alpha}W_{p_{\alpha}^{l}}-2a\Delta \mn)\rangle\nonumber\\
&-2a\alpha\langle\Delta^{s}
\left(\mn\times(\Delta \mn\times\mn)\right),
\Delta^s(\nabla_{\alpha}W_{p_{\alpha}^{l}}-2a\Delta \mn)\rangle\nonumber\\
\doteq&I_{11}+I_{12}+I_{13}+I_{14}.
\end{align}
Note that $\mh-\nabla_{\alpha}W_{p_{\alpha}^{k}}=-W_{\mn_{l}}=-2(k_3-k_2)(\mn\cdot\mcurl\mn)\mcurl\mn$, we have
\begin{eqnarray*}
I_{11}\le C_{\delta}(1+\|\nabla\mn\|_{L^{\infty}}^{2})\|\Delta^s\nabla\mn\|_{L^2}^{2}+\delta\|\Delta^{s+1}\mn\|_{L^2}^{2},
\end{eqnarray*}
and similar estimates hold for the term $I_{13}$, since $I_{13}$ can be written as the sum of a nonnegative term and a commutator term.
As to $I_{14}$, by  Lemma (\ref{lem:galiado}) and (\ref{lem:positive of h}) we have
\begin{align}\label{eq:I14}
I_{14}\le&-4a(k_1-a)\alpha\langle\nabla\Delta^s{\rm div}\mn,\nabla\Delta^s{\rm
div}\mn\rangle-4a(k_2-a)\alpha\langle\nabla\Delta^s\mcurl\mn,\nabla\Delta^s\mcurl\mn\rangle\nonumber\\
&-4a(k_3-k_2)\alpha\langle\mn\cdot\nabla_l\Delta^s\mcurl\mn,\mn\cdot\nabla_l\Delta^{s}\mcurl\mn\rangle\nonumber\\
&+C_{\delta}(\|\nabla\mn\|_{L^{\infty}}^{2}+1)\|\Delta^s\nabla\mn\|_{L^2}^{2}+\delta\|\Delta^{s+1}\mn\|_{L^2}^2\nonumber\\
\leq& C_{\delta}(\|\nabla\mn\|_{L^{\infty}}^{2}+1)\|\Delta^s\nabla\mn\|_{L^2}^{2}+\delta\|\Delta^{s+1}\mn\|_{L^2}^2.\nonumber
\end{align}
At last, we estimate $I_{12}$. Direct calculation shows that
\begin{equation}\label{nablawn}
\begin{split}
&\nabla_{\alpha}W_{p_{\alpha}^{l}}\cdot n^l\\
=&-2k_2|\nabla\mn|^2-2(k_3-k_2)(\mn\cdot\mcurl\mn)^2-2(k_1-k_2)({\rm div}\mn)^2+2(k_1-k_2)\nabla_l( n^l {\rm div}\mn).
\end{split}
\end{equation}
Thus, by Lemma \ref{lem:galiado} and (\ref{gs}) we infer that
\begin{align*}
\no I_{12}=&-2a\alpha\langle\Delta^s\nabla_{\alpha}W_{p_{\alpha}^{l}},
\Delta^{s+1}\mn\rangle+2a\alpha\langle\Delta^s((\nabla_{\alpha}W_{p_{\alpha}^{l}}\cdot n^l)\mn),\Delta^{s+1}\mn\rangle\nonumber\\
\no =&-4a^2\alpha\langle\Delta^{s+1}\mn,\Delta^{s+1}\mn\rangle-4a(k_1-a)\alpha\langle\nabla\Delta^s{\rm div}\mn,\nabla\Delta^s{\rm div}\mn\rangle\nonumber\\
\no &-4a(k_2-a)\alpha\langle\nabla\Delta^s(\nabla\times\mn),\nabla\Delta^s(\nabla\times\mn)\rangle\nonumber\\
\no &-4a(k_3-k_2)\alpha\langle\nabla_l\Delta^s(\mn\cdot(\nabla\times\mn)),\mn\cdot\nabla_l\Delta^{s}(\nabla\times\mn)\rangle\nonumber\\
\no &-4a(k_3-k_2)\alpha\langle [\nabla_l\Delta^s,\mn](\mn\cdot\nabla\times\mn),\nabla_l\Delta^{s}(\nabla\times\mn)\rangle\nonumber\\
\no &+2a\alpha\langle\Delta^s((\nabla_{\alpha}W_{p_{\alpha}^{l}}\cdot n^l)\mn),\Delta^{s+1}\mn\rangle\nonumber\\
\no \le& -4a^2\alpha\langle\Delta^{s+1}\mn,\Delta^{s+1}\mn\rangle+C_{\delta}(\|\nabla\mn\|_{L^{\infty}}^{2}+1)\|\Delta^s\nabla\mn\|_{L^2}^{2}+2\delta\|\Delta^{s+1}\mn\|_{L^2}^2,
\end{align*}
where we used Lemma (\ref{lem:galiado})-(\ref{lem:positive of h}) and  for the last term of the second equality, we have   the following observation with the help of (\ref{nablawn}):
\beno
&&\langle\Delta^s(\nabla_{l}(n^l{\rm div}\mn)\cdot n^k),\Delta^{s+1} n^k\rangle\nonumber\\
&&=\langle\Delta^s\nabla_{l}(n^l{\rm div}\mn),\Delta^{s}(n^k \Delta n^{k})\rangle-\langle\Delta^s\nabla_{l}(n^l{\rm div}\mn),[\Delta^{s},n^k] \Delta n^{k}\rangle\nonumber\\
&&\quad+\langle[\Delta^s,n^k]\nabla_{l}(n^l{\rm div}\mn),\Delta^{s+1}n^{k}\rangle.
\eeno

Similarly, we have
\begin{align*}
&2aI_2+2(k_1-a)I_2'+2(k_2-a)I_2''+2(k_3-a)I_2'''\nonumber\\
\le&C_{\delta}(\|\nabla\mn\|_{\infty}^{2}+1)\|\Delta^s\nabla\mn\|_2^2+\delta\|\Delta^{s+1}\mn\|_{2}^{2},
\end{align*}
and
\begin{eqnarray*}
|I_{3}''|+|I_4''|+|I_3'''|+|I_4'''|\le C_{\delta}(\|\nabla\mn\|_{\infty}^{2}+1)\|\Delta^s\nabla\mn\|_2^2+\delta\|\Delta^{s+1}\mn\|_{2}^{2}.
\end{eqnarray*}
Thus, the above arguments show that (\ref{eq:blow-up estimate}) is true. \endproof

\bigskip
\noindent {\bf Acknowledgments.} X. Pu is partially supported by NSFC 11471057 and Natural Science Foundation Project of CQ CSTC (cstc2014jcyjA50020). M. Wang is partially supported by  NSFC  10931001, 11371316 and Chen-Su star project by Zhejiang University. W. Wang is supported
NSFC 11301048  and "the Fundamental Research Funds for the Central Universities".

\end{document}